\documentclass[11pt]{amsart}

\numberwithin{equation}{section} 

\newcommand{\bC}{\mathbb{C}}
\newcommand{\bZ}{\mathbb{Z}}
\newcommand{\bN}{\mathbb{N}}

\newcommand{\ord}{{\mathrm{ord}}}

\newcommand{\sgn}{{\mathrm{sgn}}}

\newcommand{\Res}{\operatornamewithlimits{Res}}
\newcommand{\cN}{\mathcal{N}}

\newcommand{\cD}{\mathcal{D}}
\newcommand{\cR}{\mathcal{R}}
\newcommand{\cK}{\mathcal{K}}
\newcommand{\cQ}{\mathcal{Q}}

\newcommand{\bsnu}{\boldsymbol\nu}

\newcommand{\fX}{\mathfrak{X}}
\newcommand{\fg}{\mathfrak{g}}
\newcommand{\fh}{\mathfrak{h}}

\newcommand{\bQ}{\mathbf{ Q}}
\newcommand{\bq}{\mathbf{ q}}
\newcommand{\bK}{\mathbf{ K}}
\newcommand{\hf}{\mathbf{f}}

\newcommand{\cW}{\mathcal{W}}

\newcommand{\sa}{^{(a)}}
\newcommand{\am}{_m^{(a)}}

\newcommand{\bk}{_k^{(b)}}

\newcommand{\cj}{_j^{(c)}}

\newcommand{\aL}{^{(a)}_{L}}
\newcommand{\aLp}{^{(a)}_{L+1}}
\newcommand{\amm}{^{(a)}_{m-1}}
\newcommand{\amp}{^{(a)}_{m+1}}

\newcommand{\azero}{^{(a)}_0}

\newcommand{\hide}[1]{\null}
\newcommand{\XXX}{$X\hskip-1.2pt X\hskip-1.2pt X$}
\newcommand{\XXZ}{$X\hskip-1.2pt X\hskip-1.2pt Z$}

\usepackage{amsthm}
\newtheorem{thm}{Theorem}[section]
\newtheorem{lem}[thm]{Lemma}
\newtheorem{prop}[thm]{Proposition}
\newtheorem{cor}[thm]{Corollary}
\newtheorem{conj}[thm]{Conjecture}

\theoremstyle{definition}
\newtheorem{defn}[thm]{Definition}
\newtheorem{exmp}[thm]{Example}
\theoremstyle{remark}
\newtheorem{rem}[thm]{Remark}

\begin{document}

\title[Canonical Solutions of $Q$-systems]
{The canonical solutions of the $Q$-systems\\ and the
Kirillov-Reshetikhin conjecture\\
}

\author{Atsuo Kuniba}
\address{
Institute of Physics, University of Tokyo,
Tokyo 153-8902, Japan}
\curraddr{}
\email{atsuo@gokutan.c.u-tokyo.ac.jp}
%
\author{Tomoki Nakanishi}
\address{Graduate School of Mathematics,
    Nagoya University, Nagoya 464-8602, Japan}
\email{nakanisi@math.nagoya-u.ac.jp}
\author{Zengo Tsuboi}
\address{
Graduate School of Mathematical Sciences,
University of Tokyo,
Tokyo 153-8914, Japan}
\curraddr{}
\email{tsuboi@gokutan.c.u-tokyo.ac.jp}

\begin{abstract}
We study a class of
systems of functional equations
closely related to 
various kinds of integrable
 statistical and quantum
mechanical models.
We call them the finite and infinite $Q$-systems
according to the number of functions and
equations.
The finite $Q$-systems appear as
the thermal equilibrium conditions 
(the Sutherland-Wu equation)
for certain statistical mechanical systems.
Some infinite $Q$-systems appear as
the relations of the normalized characters of
the KR modules of the Yangians
and the quantum affine algebras.
We give two types of power series formulae
for the unique  solution
(resp.\ the unique canonical solution)
for a  finite (resp.\ infinite) $Q$-system.
As an application, we reformulate the Kirillov-Reshetikhin
conjecture on the multiplicities formula of
the KR modules
in terms of the canonical solutions of $Q$-systems.
\end{abstract}

\maketitle

\section{Introduction}
\label{sect:intr}
In the series of  works
\cite{K1,K2,KR}, Kirillov and Reshetikhin
studied the formal counting problem
 (the {\em formal completeness}) of
the Bethe vectors of the \XXX-type
integrable spin chains,
and they empirically reached
a remarkable conjectural formula on 
the characters of a certain family of
finite-dimensional modules
of the Yangian $Y(\fg)$.
Let us formulate it in the following way.

\begin{conj}
\label{conj:KR1}
Let $\fg$ be a complex simple Lie algebra
of rank $n$.
We set
$y=(y_a)_{a=1}^n$,
 $y_a=e^{-\alpha_a}$
for the simple roots $\alpha_a$ of $\fg$.
Let  $\cQ\am(y)$  be the normalized 
$\fg$-character of the KR module
$W\am(u)$ $(a=1,\dots,n$; $m=1,2,\dots$;
$u\in\bC)$ of the Yangian $Y(\fg)$;
and $\cQ^{\nu}(y)
:=\prod_{(a,m)} (\cQ\am(y))^{\nu\am}$.
Then, the formula
\begin{gather}
\label{eq:krc1}
\cQ^\nu(y){\displaystyle
\prod_{\alpha \in \Delta_+}(1-e^{-\alpha})}
=
\sum_{N=(N\am)}
\prod_{(a,m)}
\genfrac{(}{)}{0pt}{}
{ P\am(\nu,N)+N\am}{N\am}
(y_a)^{mN\am},
\\
P\am(\nu,N)=\sum_{k=1}^\infty \nu^{(a)}_k
\min(k,m)-\sum_{(b,k)}
N\bk d_ a A_{ab}\min(\frac{m}{d_b},\frac{k}{d_a})
\end{gather}
holds.
Here, $A=(A_{ab})$ is the Cartan matrix of $\fg$,
$d_a$ are coprime positive integers such that 
$(d_a A_{ab})$ is symmetric,
$\Delta_+$ is the set of all the positive
roots of $\fg$,
and $\binom{a}{b}=
\varGamma(a+1)/\varGamma(a-b+1)\varGamma(b+1)$.
\end{conj}

\begin{rem}
Due to the Weyl character formula,
the series in the RHS of (\ref{eq:krc1}) should be
a {\em polynomial\/} of $y$,
and its coefficients are identified with
the multiplicities of the $\fg$-irreducible
components of the tensor product
$\bigotimes_{(a,m)}W\am(u\am)^{\otimes \nu\am}$,
where $u\am$ are arbitrary.
\end{rem}

\begin{rem}
\label{rem:binom1}
There are actually two versions of
Conjecture \ref{conj:KR1}.
The above one is the version in \cite{HKOTY}
which followed \cite{K1,K2}.
In the version in \cite{KR},
the binomial coefficients $\binom{a}{b}$ are set to be 0 if $a<b$;
furthermore, the equality is claimed,
not for the entire series
in the both hand sides of (\ref{eq:krc1}),
but only for their coefficients of the powers $y^M$
``in the fundamental Weyl chamber'';
namely, $M=(M_a)_{a=1}^n$ satisfies
\begin{align}
\sum_{(a,m)}\nu\am m\Lambda_a-\sum_{a=1}^n
M_a \alpha_a \in P_+,
\end{align}
where $\Lambda_a$ are the fundamental weights
and $P_+$ is the set of the dominant
integral weights of $\fg$.
So far, it is not proved that the two conjectures
are equivalent.
The both conjectures are naturally
translated
into the ones  for the untwisted
quantum affine algebras, which are extendable
to the twisted quantum affine algebras \cite{HKOTT}.
In this paper,
we refer all these conjectures as the Kirillov-Reshetikhin
conjecture.
More comments 
and the current status of the conjecture
will be given in Section  \ref{subsec:rel1}.
\end{rem}

In \cite{KR,K3},
it was claimed that
the $\cQ\am(y)$'s satisfy
a system of equations
\begin{align}\label{eq:qsys23}
\begin{split}
(\cQ\am(y))^2
&=
\cQ\amm(y)\cQ\amp(y)\\
&\qquad\qquad
+(y_a)^m (\cQ\am(y))^2\prod_{(b,k)}
(\cQ\bk(y))^{G_{am,bk}}.
\end{split}
\end{align}
Here, $\cQ\azero(y)=1$, and $G_{am,bk}$ are the integers
defined as
\begin{align}
\label{eq:bformula2}
G_{am,bk}=
\begin{cases}
-A_{ba}
(\delta_{m,2k-1}+2\delta_{m,2k}+\delta_{m,2k+1})
&
d_b/d_a=2\\
-A_{ba}
(\delta_{m,3k-2}+2\delta_{m,3k-1}+3\delta_{m,3k}&
d_b/d_a=3\\
\qquad\qquad\qquad
+2\delta_{m,3k+1}+\delta_{m,3k+2})&\\
-A_{ab}\delta_{d_a m, d_b k}
& \text{otherwise}.
\end{cases}
\end{align}
See (\ref{eq:qsys25}) for the original form of
(\ref{eq:qsys23}) in
\cite{KR,K3}.
The relations (\ref{eq:qsys23}) and (\ref{eq:qsys25})
are  often called the {\em $Q$-system}.
The importance of the role of the
$Q$-system to the
formula (\ref{eq:krc1}) was
recognized in  \cite{K1,K2,KR},
and  more explicitly exhibited
in \cite{HKOTY,KN2}.
In this paper we proceed one step
further in this direction;
we study the equation (\ref{eq:qsys23})
in a  more general point of view,
and give a characterization of 
the special power series
solution in (\ref{eq:krc1}).
For this purpose,
 we introduce
{\em finite and infinite $Q$-systems},
where the former (resp.\ the latter) is a finite 
(resp.\ infinite) system of equations
for a finite (resp.\ infinite) family of
power series of the variable with
finite (resp.\ infinite) components.
The equation (\ref{eq:qsys23}),
which is an infinite system of equations
with the variable with finite components,
is regarded as
an infinite $Q$-system with the specialization
of the variable (a {\em specialized $Q$-system}).
We show that every finite $Q$-system has a unique
solution which has the same type of
the power series formula as (\ref{eq:krc1})
(Theorem \ref{thm:fqmain1}).
In contrast,
infinite $Q$-systems and their specializations,
 in general, admit more than one solutions.
However, every infinite $Q$-system,
or its specialization, has a unique
{\em canonical\/} solution (Theorems \ref{thm:cs4}
and \ref{thm:cs2}),
whose definition is given in Definition
\ref{def:can1}.
The formula  (\ref{eq:krc1})
turns out to be exactly the power series formula
for the canonical solution 
of (\ref{eq:qsys23}) (Theorem \ref{thm:psf5}
and Proposition \ref{prop:denom11}).
Therefore, one can rephrase Conjecture
\ref{conj:KR1} in a more intrinsic way
as follows (Conjecture \ref{conj:KR2}):
{\em The family $(\cQ\am(y))$ of the normalized
$\fg$-characters of the KR modules is characterized
as the canonical
solution of (\ref{eq:qsys23}).}
This is the main statement of the paper.

Interestingly,
the finite $Q$-systems  also appear
in other types of 
integrable statistical mechanical systems.
Namely, they appear as
the thermal equilibrium condition
(the Sutherland-Wu equation)
for the Calogero-Sutherland model
\cite{S}, as well as
the one for the ideal gas of the Haldane
exclusion statistics  \cite{W}.  
The property of the solution
of the finite $Q$-systems are studied 
in \cite{A,AI,IA} from the point of view
of the quasi-hypergeometric functions.
We expect that the study of the $Q$-system
and its variations and  extensions
will be useful for the
representation theory of the quantum groups,
and for the understanding of the nature
of the integrable models as well.

We would like to thank 
V.\ Chari,
G.\ Hatayama,
A.\ N.\ Kirillov,
M.\ Noumi,
M.\ Okado,
T.\ Takagi,
and Y.\ Yamada
for
very useful discussions.
We especially thank K.\ Aomoto
for the discussion where we recognize the
very close relation between the present work
and his work,
and also for pointing out the reference \cite{G} to us.

\section{Finite $Q$-systems}
\label{sect:finiteQ}

A considerable part of the  results 
in this section 
can be found in the work by
Aomoto and Iguchi  \cite{A,IA}.
We present here a more direct approach.
More detailed remarks will be given in
 Section \ref{subsec:hal1}.

\subsection{Finite $Q$-systems}
Throughout Section \ref{sect:finiteQ}, let $H$ denote a finite
index set.
Let $w=(w_i)_{i\in H}$ and $v=(v_i)_{i\in H}$
be complex multivariables,
and let
$G=(G_{ij})_{i,j\in H}$
be a given complex square matrix of size $|H|$.
We  consider a holomorphic map
$\cD\rightarrow \bC^H$,
$v\mapsto w(v)$ with
\begin{align}
\label{eq:wv1}
w_i(v)&=v_i\prod_{j\in H} (1-v_j)^{-G_{ij}},
\end{align}
where $\cD$ is some neighborhood of $v=0$
 in $\bC^H$.
The Jacobian $(\partial w / \partial v)(v) $ is 1
at $v = 0$, so that
the map $w(v)$
is bijective around $v=w=0$.
Let $v(w)$ be the inverse map around $v=w=0$.
Inverting (\ref{eq:wv1}),
we obtain the following functional equation
for $v_i(w)$'s:
\begin{align}\label{eq:vw1}
v_i(w)
= w_i\prod_{j\in H}(1-v_j(w))^{G_{ij}}.
\end{align}
By introducing new functions
\begin{align}\label{eq:qv1}
Q_i(w)=1-v_i(w),
\end{align}
the equation
(\ref{eq:vw1}) is written as
\begin{align}\label{eq:qsys1}
Q_i(w)+w_i \prod_{j\in H}(Q_j(w))^{G_{ij}}=1.
\end{align}

\par
{}From now on, we regard 
(\ref{eq:qsys1}) as a system of
equations for a family $(\cQ_i(w))_{i\in H}$ of
power series  of
$w=(w_i)_{i\in H}$
with the unit constant terms
({\em i.e.,} the constant terms are $1$).
Here, for any power series $f(w)$ with the unit constant term
and any complex number $\alpha$,
we mean by  $(f(w))^\alpha\in \bC[[w]]$
the $\alpha\/$th power of $f(w)$ with the unit constant term.
We can easily reverse 
the procedure from (\ref{eq:wv1}) to (\ref{eq:qsys1}),
and we have
\par
\begin{prop}\label{prop:qsys1}
The power series expansion of $Q_i(w)$ in (\ref{eq:qv1})
gives the unique family 
$(Q_i(w))_{i\in H}$ of power series of
$w$ with the unit constant terms
which satisfies (\ref{eq:qsys1}).
\end{prop}

\begin{defn}
The following system of equations for
a family $(Q_i(w))_{i\in H}$ of 
 power series of $w$ with the unit constant terms is called
a {\em (finite) $Q$-system}:
For each $i\in H$,
\begin{align}\label{eq:qsys2}
\prod_{j\in H}
(Q_j(w))^{D_{ij}}+w_i \prod_{j\in H}(Q_j(w))^{G_{ij}}=1,
\end{align}
where $D=(D_{ij})_{i,j\in H}$
 and $G=(G_{ij})_{i,j\in H}$ are arbitrary
 complex matrices with $\det D\neq 0$.
The equation (\ref{eq:qsys1}), which is the special case
of (\ref{eq:qsys2}) with $D=I$ ($I$: the identity matrix),
is called a {\em standard\/} $Q$-system.
\end{defn}

It is easy to see that
there is a one-to-one correspondence between
the solutions of the $Q$-system (\ref{eq:qsys2})
and the solutions of the standard $Q$-system
\begin{align}\label{eq:qsys3}
Q'_i(w)+w_i \prod_{j\in H}(Q'_j(w))^{G'_{ij}}=1,
\quad G'=GD^{-1},
\end{align}
where the correspondence is given by
\begin{align}\label{eq:qd1}
Q'_i(w)&=\prod_{j\in H}
(Q_j(w))^{D_{ij}},\\
\label{eq:qd8}
Q_i(w)&=\prod_{j\in H}
(Q'_j(w))^{(D^{-1})_{ij}}.
\end{align}

Therefore, from Proposition \ref{prop:qsys1},
 we immediately have
\begin{thm}
\label{thm:qunique1}
There exists a unique solution
of the $Q$-system
(\ref{eq:qsys2}), which is given by
(\ref{eq:qd8}),
where $(Q'_i(w))_{i\in H}$ is the unique solution of
the standard $Q$-system (\ref{eq:qsys3}).
\end{thm}

\subsection{Power series formulae}
\label{subsec:psf1}

\label{subsec:pse1}
In what follows,
we use the binomial coefficient in the following sense:
For $a\in \bC$ and $b\in \bZ_{\geq 0}$,
\begin{align}
\label{eq:binom1}
\binom{a}{b}
=\frac{\varGamma(a+1)}{\varGamma(a-b+1)\varGamma(b+1)},
\end{align}
where the RHS means the limit value for the singularities.
We set $\cN:=(\bZ_{\geq 0})^H$.
For  $D$, $G$ in (\ref{eq:qsys2})
and $\nu=(\nu_i)_{i\in H} \in \bC^{H}$,
we define two power series of $w$,
\begin{align}
\label{eq:cKdef1}
K^\nu_{D,G}(w) &= \sum_{N\in \cN} K(D,G;\nu,N) w^N,
\quad w^N=\prod_{i\in H}w_i^{N_i},\\
\label{eq:cRdef1}
R^\nu_{D,G}(w) &= \sum_{N\in \cN} R(D,G;\nu,N)w^N
\end{align}
with the coefficients
\begin{align}
\label{eq:kdef1}
K(D,G;\nu,N) &= \prod_{i\in H(N)}
\binom{P_i + N_i}{N_i},\\
\label{eq:rdef1}
R(D,G;\nu,N) &=
\Bigl(\det_{H(N)} F_{ij}\Bigr)
\prod_{i\in H(N)} \frac{1}{N_i}
\binom{P_i+ N_i - 1}{N_i- 1},
\end{align}
where we set $H(N)=\{\, i\in H\mid N_i\neq 0\, \}$
for each $N\in \cN$,
\begin{align}\label{eq:padef1}
P_i&=P_i(D,G;\nu,N):=-
\sum_{j\in H}\nu_j(D^{-1})_{ji}
-\sum_{j\in H}N_j(GD^{-1})_{ji},\\
\label{eq:fdef1}
F_{ij}&=
F_{ij}(D,G;\nu,N):=
\delta_{ij}P_j+(GD^{-1})_{ij}N_j,
\end{align}
and $\det_{H(N)}$ is a shorthand notation for $\det_{i,j\in H(N)}$.
In (\ref{eq:kdef1}) and (\ref{eq:rdef1}),
 $\det_\emptyset$ and $\prod_\emptyset$
mean 1; namely, 
$K^\nu_{D,G}(w)$ and 
$R^\nu_{D,G}(w)$ are power series
with the unit constant terms.
It is easy to check that
the both series converge for
 $|w_i|<|\gamma_i^{\gamma_i}/(\gamma_{i}+1)^{\gamma_{i}+1}|$,
where $\gamma_i=-(GD^{-1})_{ii}$ and
$z^z=\exp(z\log z)$ with the principal
 branch $-\pi < \mathrm{Im} (\log z) \leq \pi$ chosen.

\par
Now we state our main results in this section.

\begin{thm}[Power series formulae]\label{thm:fqmain1}
 Let $(Q_i(w))_{i\in H}$  be the unique
solution of (\ref{eq:qsys2}).
For $\nu\in \bC^{H}$, let
$Q^\nu_{D,G}(w):=\prod_{i\in H}
(Q_i(w))^{\nu_i}$.
Then,
\begin{align}
\label{eq:qkr01}
Q^\nu_{D,G}(w)&=K^\nu_{D,G}(w)/
K^0_{D,G}(w),\\
\label{eq:qkr02}
Q^\nu_{D,G}(w)&
=R^\nu_{D,G}(w).
\end{align}
\end{thm}

The power series formulae for $Q_i(w)$ are
obtained as special cases of (\ref{eq:qkr01})
and (\ref{eq:qkr02}) by setting $\nu=(\nu_j)_{j\in H}$
as $\nu_j=\delta_{ij}$.

One may recognize that
the first formula (\ref{eq:qkr01}) is analogous to the
formula  (\ref{eq:krc1}),
where the denominator $K^0_{D,G}(w)$ in 
(\ref{eq:qkr01})
corresponds to the Weyl denominator in the LHS
of (\ref{eq:krc1}).
As mentioned in Section \ref{sect:intr},
the formula (\ref{eq:krc1}) is interpreted as
the formal completeness of the {\em \XXX-type\/}
Bethe vectors.
In the same sense,
the second formula (\ref{eq:qkr02})
is analogous to the  formal
completeness of the {\em \XXZ-type\/} Bethe vectors
in \cite{KN1,KN2}.
See Section \ref{subsec:hal1} for more remarks.

\begin{exmp}
Let $|H|=1$.
Then, (\ref{eq:qsys2}) is an equation
for a single power series $Q(w)$,
\begin{align}
\label{eq:lamb1}
(Q(w))^D+w(Q(w))^G =1,
\end{align}
where $D\neq 0$ and $G$ are complex numbers,
and the series (\ref{eq:cRdef1}) reads as
\begin{align}
\label{eq:lamb2}
R^\nu_{D,G}(w)
=\frac{\nu}{D}
\sum_{N=0}^\infty
\frac{\varGamma((\nu+NG)/D)(-w)^N}
{\varGamma((\nu+NG)/D-N+1)N!}.
\end{align}
The equation (\ref{eq:lamb1}) 
and the power series formula (\ref{eq:lamb2})
are well known and
have a very long history
since Lambert  (e.g.\ \cite[pp.\ 306--307]{B}).
\end{exmp}

\begin{exmp}
Consider the case $G=O$ in (\ref{eq:qsys2}),
\begin{align}
\prod_{j\in H}
(Q_j(w))^{D_{ij}}+w_i =1.
\end{align}
This is easily solved as
\begin{align}
\label{eq:qe1}
Q_i(w)=
\prod_{j\in H}
(1-w_j)^{(D^{-1})_{ij}},
\end{align}
and, therefore,
\begin{align}
\label{eq:qe2}
Q^\nu_{D,O}(w)=
\prod_{i\in H}
(1-w_i)^{\sum_{j\in H}\nu_j (D^{-1})_{ji}}
=
\prod_{i\in H}
(1-w_i)^{-P_i(D,O;\nu,N)},
\end{align}
where $N\in \cN$ is arbitrary.
Using the binomial theorem
\begin{equation}\label{eq:exp}
(1-x)^{-\beta-1} = \sum_{N=0}^\infty \binom{\beta+N}{N}x^N,
\end{equation}
one can directly check that
\begin{align}
Q^\nu_{D,O}(w)&=\sum_{N\in \cN}
\prod_{i\in H(N)}\binom{P_i-1+N_i}{N_i}w_i^{N_i}
=R^{\nu}_{D,O}(w),\\
Q^\nu_{D,O}(w)&=
\frac{\prod_{j\in H}(1-w_i)^{\sum_{j\in H}\nu_j 
(D^{-1})_{ji}-1}}
{\prod_{j\in H}(1-w_i)^{-1}}
=
\frac{K^\nu_{D,O}(w)}
{K^{0}_{D,O}(w)}.
\end{align}
\end{exmp}

\subsection{Proof of Theorem \ref{thm:fqmain1}
and basic formulae}
\label{subsec:pf1}

Theorem \ref{thm:fqmain1} is regarded as
a particularly
nice example of the multivariable Lagrange
inversion formula (e.g.\ \cite{G})
where all the explicit calculations
can be carried through.
Here, we present the most direct calculation
based on
the multivariable residue formula
(the {\em Jacobi formula\/} in \cite[Theorem 3]{G}).

We first remark that
\begin{lem}
\label{lem:kkd1}
Let $G'=GD^{-1}$. For each $\nu\in \bC^{H}$,
let $\nu'\in \bC^{H}$
with
$\nu'_i=\sum_{j\in H}\nu_j (D^{-1})_{ji}$.
Then,
\begin{align}
\label{eq:qq6}
Q^{\nu}_{D,G}(w)&=
Q^{\nu'}_{I,G'}(w),\\
\label{eq:kk6}
K^{\nu}_{D,G}(w)&=
K^{\nu'}_{I,G'}(w),
\quad
R^{\nu}_{D,G}(w)=
R^{\nu'}_{I,G'}(w).
\end{align}
\end{lem}
\begin{proof}
The equality
(\ref{eq:qq6}) is due to Theorem  \ref{thm:qunique1}.
The ones (\ref{eq:kk6})
 follow from the fact $P_i(D,G;\nu,N)=P_i(I,G';\nu',N)$.
\end{proof}

By Lemma \ref{lem:kkd1},
we have only to prove 
Theorem \ref{thm:fqmain1}
for the standard case $D=I$.
Recall that (Proposition \ref{prop:qsys1})
 $Q^\nu_{I,G}(w)=
\prod_{i\in H}(1-v_i(w))^{\nu_i}$,
where $v=v(w)$ is the inverse map of (\ref{eq:wv1}).
Thus, Theorem \ref{thm:fqmain1}
follows from
\begin{prop}[Basic formulae]
\label{prop:pse2}
Let $v=v(w)$ be the inverse map of (\ref{eq:wv1}).
Then, the power series expansions
\begin{align}
\label{eq:pse1}
\det_H \Bigl(\frac{w_j}{v_i}
\frac{\partial v_i}{\partial w_j}(w)\Bigr)
\prod_{i\in H}
(1-v_i(w))^{\nu_i-1}
&=K^\nu_{I,G}(w),\\
\label{eq:vr1}
\prod_{i\in H}(1-v_i(w))^{\nu_i}
&=
R^\nu_{I,G}(w)
\end{align}
hold around $w=0$.
\end{prop}

\begin{proof}
{\em The first formula (\ref{eq:pse1}).} 
\noindent
%
We evaluate the coefficient for 
$w^N$ in the LHS of (\ref{eq:pse1}) as follows:
\begin{align*}
&\,
\Res_{w=0} 
\frac{\partial v}{\partial w}(w)
\prod_{i\in H}
\Bigl\{ (1-v_i(w))^{\nu_i-1} (v_i(w))^{-1}
(w_i)^{1-N_i-1} \Bigr\}dw\\
=&\, \Res_{v=0} 
\prod_{i\in H}
\Bigl\{ (1-v_i)^{\nu_i-1}(v_i)^{-1}
\Bigl(v_i\prod_{j\in H}(1-v_j)^{-G_{ij}}\Bigr)^{-N_i}
 \Bigr\}dv\\
=&\, \Res_{v=0} 
\prod_{i\in H}
\Bigl\{ (1-v_i)^{-P_i(I,G;\nu,N)-1}
(v_i)^{-N_i-1}
 \Bigr\}dv\\
= &\, \prod_{i\in H}
\binom{P_i(I,G;\nu,N) + N_i}{N_i}
=K(I,G;\nu,N),
\end{align*}
where we used (\ref{eq:exp}) to get the last line.
Thus,  (\ref{eq:pse1}) is proved.

\par
{\em The second formula (\ref{eq:vr1}).}
By a simple calculation, we have
\begin{align}
\begin{split}
\det_H \Bigl(\frac{v_j}{w_i}
\frac{\partial w_i}{\partial v_j}(v)\Bigr)
\prod_{i\in H}(1-v_i)
&=\det_{ H}
\Bigl(\delta_{ij}+(-\delta_{ij}+G_{ij})v_i\Bigr)
\\
\label{eq:k02}
&= \sum_{J\subset H}
d_J \prod_{i\in J}v_i,
\end{split}
\end{align}
where $d_J:= \det_J (-\delta_{ij}+G_{ij})$,
and the sum is taken over
all the subsets $J$ of $H$.
Therefore,
the LHS of (\ref{eq:vr1})
is written as
\begin{align}
\label{eq:qig1}
\det_H \Bigl(\frac{w_j}{v_i}
\frac{\partial v_i}{\partial w_j}(w)\Bigr)
 \sum_{J\subset H}
 d_J 
\prod_{i\in H}
\Bigl\{
(1-v_i(w))^{\nu_i-1}
v_i(w)^{\theta(i\in J)}\Bigr\}.
\end{align}
By a similar residue calculation as above,
the  coefficient for $w^N$ of (\ref{eq:qig1})  is evaluated
as ($\theta(\text{\rm true}) = 1$ and
$\theta(\text{\rm false}) = 0$)
\begin{align*}
&\, 
\sum_{J\subset H}
d_J
\Res_{v=0} 
\prod_{i\in H}
\Bigl\{ (1-v_i)^{-P_i(I,G;\nu,N)-1}
(v_i)^{-N_i+\theta(i\in J)-1}
 \Bigr\}dv\\
= &\, 
\sum_{J\subset H(N)}d_J
\prod_{i\in H(N)}
\binom{P_i(I,G;\nu,N) + N_i-\theta(i\in J)}{N_i-\theta(i\in J)}\\
=&\,
\biggl(
\sum_{J \subset H(N)}
d_J
\prod_{i \in  J} N_i
\prod_{i \in H(N)\setminus J}(P_i+N_i)
\biggr)
\prod_{i\in H(N)} \frac{1}{N_i}
\binom{P_i+ N_i - 1}{N_i- 1}\\
=&\,
\det_{H(N)} \Bigl(\delta_{ij}(P_j+N_j)+(-\delta_{ij}+G_{ij})N_j
\Bigr)
\prod_{i\in H(N)} \frac{1}{N_i}
\binom{P_i+ N_i - 1}{N_i- 1}\\
=&\, R(I,G;\nu,N).
\end{align*}
Thus, (\ref{eq:vr1}) is proved.
\end{proof}

This completes the proof of Theorem \ref{thm:fqmain1}.

\begin{exmp}
\label{exmp:tri1}
We say that the map $w(v)$ in (\ref{eq:wv1}) is
{\em lower-triangular\/}
if the matrix
$G_{ij}$ is strictly lower-triangular
with respect to 
a certain total order $\prec$
in $H$ ({\em i.e.,} $G_{ij}=0$ for $i\preceq j$).
Let $w(v)$ be a lower-triangular map.
Then,
\begin{align}
\label{eq:dvw1}
\det_H \Bigl(\frac{v_j}{w_i}
\frac{\partial w_i}{\partial v_j}(v)\Bigr)
=\det_H\Bigl(\delta_{ij}+\frac{G_{ij}v_j}{1-v_j}\Bigr)
=1.
\end{align}
Thus, the formula (\ref{eq:pse1})
is simplified as
\begin{align}
\label{eq:pse9}
\prod_{i\in H}
(1-v_i(w))^{\nu_i-1}
&=K^\nu_{I,G}(w).
\end{align}
This type of formulae has appeared in \cite{K1,K2,HKOTY}.
\end{exmp}

Let us isolate the case $\nu=0$
{}from (\ref{eq:pse1}),
together with the formula (\ref{eq:k02}),
for the later use:

\begin{cor}[Denominator formulae]
\begin{align}
\label{eq:pse2}
K^0_{I,G}(w)&=
\det_H \Bigl(\frac{w_j}{v_i}
\frac{\partial v_i}{\partial w_j}(w)\Bigr)
\prod_{i\in H}
 (1-v_i(w))^{-1},\\
\label{eq:den1}
K^0_{I,G}(w)&=
\Bigl\{
\det_{ H}
\Bigl(\delta_{ij}(1-v_i(w))+G_{ij}v_i(w)\Bigr)
\Bigl\}^{-1}.
\end{align}
\end{cor}

{}From (\ref{eq:den1}) and 
the first formula of Theorem \ref{thm:fqmain1},
we obtain
\begin{cor}
\label{cor:qk3}
\begin{align}
\label{eq:qk1}
Q^\nu_{I,G}(w)
=&\,\sum_{J\subset H}
g_J K^{\nu+\delta_J}_{I,G}(w),\\
g_J:=&\,
\sum_{\genfrac{}{}{0pt}{1}{J'\subset H}{|J'|=|J|}}
\sgn \binom{J\overline{J}}{J'\overline{J'}}
\det_{i\in J,j\in J'}
\bigl(\delta_{ij}-G_{ij}\bigr)
\det_{i\in \overline{J},j\in \overline{J'}}
G_{ij},
\end{align}
where $\delta_J=(\theta_i)_{i\in H}$,
$\theta_i=1$ if $i\in J$ and $0$ otherwise,
and $\overline{J}=H\setminus J$.
\end{cor}

{}From  Corollary \ref{cor:qk3}, one can easily
reproduce the second formula
of Theorem \ref{thm:fqmain1}.
We leave it as an exercise for the reader.

\subsection{Remarks on related works}
\label{subsec:hal1}

{\em i) The formal completeness of
the  Bethe vectors.}
In \cite{K1,K2,HKOTY,KN1,KN2,KNT},
the formal completeness
of the \XXX/\XXZ-type Bethe vectors
are studied.
In the course of their analysis,
several power series formulae
in this section
appeared
in specialized/implicit forms.
For example,   Lemma 1 in \cite{K1}
is a special case of  (\ref{eq:pse9}),
Theorem 4.7 in \cite{KN2} is a special case
of Proposition \ref{prop:pse2},
{\em etc}.
{}From the current point of view, however,
the relation between
these power series
formulae and the underlying  {\em finite\/} $Q$-systems
was not clearly recognized therein.
As a result,
these power series formulae and the {\em infinite\/} $Q$-systems
were somewhat abruptly combined in the limiting procedure
to obtain the power series formula for the {\em infinite\/}
$Q$-systems.
We are going to straighten out this logical entangle,
and make the logical structure more transparent by
Theorem \ref{thm:fqmain1}
and the forthcoming Theorems \ref{thm:qkr3},
\ref{thm:psf5}, Proposition \ref{prop:denom11},
and Conjecture \ref{conj:KR2}.

\par
{\em ii) The ideal gas with the Haldane
statistics and the Sutherland-Wu equation.}
The series $K^\nu_{D,G}(w)$ has an interpretation
of the grand partition function of the ideal
gas with the Haldane  exclusion statistics \cite{W}.
The finite $Q$-system appeared
in \cite{W} as
the thermal equilibrium condition
for the distribution functions
of the same system.
See also \cite{IA} for another 
interpretation.
The one variable case (\ref{eq:lamb1})
also appeared in \cite{S}
as the thermal equilibrium condition
for the distribution function
of the Calogero-Sutherland model.
As an application of our second formula
in Theorem \ref{thm:fqmain1},
we can quickly reproduce
 the ``cluster expansion
formula'' in \cite[Eq.\ (129)]{I},
which was originally calculated by the Lagrange
inversion formula, as follows:
\begin{align}
\begin{split}
 &\log Q_i(w)
=
\Bigl[
\frac{\partial}{\partial \nu_i}
R^\nu_{I,G}(w)\Bigr]_{\nu=0}\\
=&\
\sum_{N\in \cN}
\det_{\genfrac{}{}{0pt}{1}{H(N)}{j,k\neq i}} F_{jk}(I,G;0,N)
\prod_{j\in H(N)}
\frac{1}{N_j}
\binom{P_j(I,G;0,N)+N_j-1}{N_j-1}
w^N,
\end{split}
\end{align}
where $\{Q_i(w)\}_{i\in H}$ is
the solution of (\ref{eq:qsys1}).
The Sutherland-Wu equation also
plays an important role
for the conformal field theory spectra.
(See \cite{BS} and the references therein.)

\par
{\em iii) Quasi-hypergeometric functions.}
The series $K^\nu_{D,G}(w)$
is a special example of the 
quasi-hypergeometric
functions by Aomoto and Iguchi \cite{AI}; 
when $G'_{ij}$ are all integers,
it reduces to a 
 general hypergeometric
function of Barnes-Mellin type.
A quasi-hypergeometric function satisfies
a system of fractional differential equations
and a system of difference-differential equations
\cite{AI}.
It also admits an integral representation \cite{A}.
In particular,
the integral representation for $K^\nu_{I,G}(w)$
reduces to a simple form 
(\cite[Eq.\ (2.30)]{A}, \cite[Eq.\ (89)]{IA});
in our notation,
\begin{align}
\label{eq:ao2}
K^\nu_{I,G}(w)
&=
\frac{1}{(2\pi \sqrt{-1})^{|H|}}
\int
\Bigl\{
 \prod_{i\in H}t_i^{\nu_i-1}
f_i(w,t)^{-1}\Bigr\}dt,\\
f_i(w,t)&:=t_i-1+w_i\prod_{j\in H}t_j^{G_{ij}},
\end{align}
where the integration is along a circle
around $t_i=1$ starting from $t_i=0$ for each $t_i$.
We see that $f_i(w,t)=0$ is the standard $Q$-system
(\ref{eq:qsys1}).
The integral (\ref{eq:ao2}) is easily
evaluated by the
Cauchy theorem as \cite[eq.\ (2.32)]{A}
\begin{align}
\label{eq:ao1}
K^\nu_{I,G}(w)
=
Q^\nu_{I,G}(w)/\det_{ H}
(\delta_{ij}Q_i(w)
 + G_{ij}(1-Q_j(w))),
\end{align}
where $\{ Q_i(w)\}_{i\in H}$ is the solution
of (\ref{eq:qsys1}).
The formula (\ref{eq:ao1}) reproduces a version
of the Lagrange inversion formula
(the Good formula \cite[Theorem 2]{G}),
and it is equivalent to
the formulae (\ref{eq:qkr01}),
 (\ref{eq:k02}), and  (\ref{eq:pse2}).

\section{Infinite $Q$-systems}
\label{sec:infiniteQ}

\subsection{Infinite $Q$-systems}
Throughout Section \ref{sec:infiniteQ},
let $H$ be a countable infinite index set.
We fix an increasing sequence of
{\em finite\/} subsets  of $H$,
$H_1\subset H_2\subset \cdots \subset H$
such that $\varinjlim H_L=H$.
The result below does not depend on the choice
of the sequence $\{H_L\}_{L=1}^\infty$.
A natural choice is  $H=\bN$ and $H_L=\{\, 
1,\dots, L\, \}$. 
However, we introduce this generality
to accommodate the situation we encounter in Section
\ref{sec:spec} (cf.\ (\ref{eq:hinf1})).

Let $w=(w_i)_{i\in H}$ be a multivariable
with infinitely many components.
For each $L\in \bN$,
let $w_L=(w_i)_{i\in H_L}$ be the
submultivariable of $w$.
The field $\bC[[w_L]]$
of the power series
 of $w_L$ over $\bC$
is equipped with the standard
$\fX_L$-adic topology,
where $\fX_L$ is the ideal of $\bC[[w_L]]$
generated by $w_i$'s ($i\in H_L$).
For $L<L'$, there is a natural projection
$p_{LL'}:\bC[[w_{L'}]]\rightarrow \bC[[w_L]]$
such that $p_{LL'}(w_i)=w_i$  if $i\in H_L$
and 0 if $i\in H_{L'}\setminus H_L$.
A {\em power series $f(w)$ of $w$} is an element
of the projective limit
 $\bC[[w]]=\varprojlim \bC[[w_L]]$ 
of the projective system
\begin{align}
\label{eq:pl1}
\bC[[w_1]]
\leftarrow \bC[[w_2]]\leftarrow \bC[[w_{3}]]
\leftarrow\cdots
\end{align}
with the induced topology.
Let $p_L$ be the canonical projection $p_L:\bC[[w]]\rightarrow
\bC[[w_L]]$,
and $f_L(w_L)$ be the $L$th projection image of $f(w)\in \bC[[w]]$;
 namely, 
$f_L(w_L)=p_L(f(w))$ and
$f(w)=(f_L(w_L))_{L=1}^\infty$.

Here are some basic properties of power series which we use below:
\par
(i) We also present a power series $f(w)$ as a formal sum
\begin{gather}
f(w)=\sum_{N\in \cN} a_Nw^N,
\quad a_N\in \bC,
\label{eq:f1}\\
\cN=\{\, N=(N_i)_{i\in H}
\mid N_i\in \bZ_{\geq 0},
\, \text{all but finitely many $N_i$ are zero}\,\},
\end{gather}
(the definition of $\cN$  is reset here
for the infinite index set $H$)
whose $L$th projection image is
\begin{gather}
\label{eq:fl1}
f_L(w_L)=\sum_{N\in \cN_L} a_Nw^N,\\
\label{eq:nldef1}
\mathcal{N}_L=
\{\, N\in \cN \mid
N_i = 0\ \text{for}\ i\notin H_L\,\}.
\end{gather}

\par
(ii)
For any  power series $f(w)$ with the unit constant term
and any complex number $\alpha$,
the $\alpha\/$th power $(f(w))^\alpha
:=((f_L(w_L))^\alpha)_{L=1}^\infty \in \bC[[w]]$ is
uniquely defined and has the unit constant term again.

\par
(iii)
Let  $f_i(w)$ ($i\in H$) be a family of power series
and $f_{i,L}(w_L)$ be their $L$th projections.
If their  infinite product exists
irrespective of the order of the product,
we write it as $\prod_{i\in H} f_{i}(w)$.
$\prod_{i\in H} f_{i}(w)$ exists
if and only if $\prod_{i\in H} f_{i,L}(w_L)$ exists for each $L$;
furthermore,
if they exist, the latter is the $L$th projection of the former.

\begin{defn}\label{def:qsys1}
The following system of equations for
a family $(Q_i(w))_{i\in H}$ of 
 power series of $w$ with the unit constant terms
is called
an {\em (infinite) $Q$-system}:
For each $i\in H$,
\begin{align}\label{eq:qsys8}
\prod_{j\in H}
(Q_j(w))^{D_{ij}}+w_i \prod_{j\in H}(Q_j(w))^{G_{ij}}=1.
\end{align}
Here, $D=(D_{ij})_{i,j\in H}$
 and $G=(G_{ij})_{i,j\in H}$ are arbitrary
infinite-size complex matrices
satisfying the following two conditions:
\begin{itemize}
\item[(D)] The matrix $D$ is invertible, {\em i.e.},
there exists a matrix $D^{-1}$ such that
$DD^{-1}=D^{-1}D=I$.
\item[(G')] The matrix product $G'=GD^{-1}$ is well-defined.
\end{itemize}
When $D=I$, the equation (\ref{eq:qsys8})
 is called a {\em standard\/} $Q$-system.
\end{defn}

\begin{rem}
\label{rem:dg1}
The condition (G') is rephrased as
``for each  $i$ and $k$, 
all but finitely many 
$G_{ij}(D^{-1})_{jk}$ ($j\in H$) are zero''.
Similarly, the condition (D) implies that,
for each  $i$ and $k$, 
all but finitely many 
$D_{ij}(D^{-1})_{jk}$, $(D^{-1})_{ij}D_{jk}$
 ($j\in H$) are zero.
For the standard case,
(D) is trivially satisfied,
and (G') is satisfied for any complex matrix $G$.
\end{rem}

Unlike the finite $Q$-systems,
the uniqueness of the solution does not
hold for the infinite $Q$-systems, in general.
For instance, the following example
admits infinitely many solutions.

\begin{exmp}
\label{exmp:inf1}
Let $H=\bZ$, and
consider a $Q$-system,
\begin{align}
\label{eq:qi2}
\frac{Q_{i-1}(w)Q_{i+1}(w)}
{(Q_i(w))^2} + w_i =1,
\end{align}
where $Q_0(w)=1$. 
This can be easily solved as
\begin{align}
\label{eq:qi1}
Q_i(w)= (Q_1(w))^i\prod_{j=1}^{i-1}(1-w_j)^{i-j},
\end{align}
where $Q_1(w)$ is an arbitrary series of $w$
with the unit constant term.
\end{exmp}

\subsection{Canonical solution}

\subsubsection{Solution of standard $Q$-system}

First, we consider the standard case
\begin{align}\label{eq:qsys5}
Q_i(w)+w_i \prod_{j\in H}(Q_j(w))^{G_{ij}}=1.
\end{align}
Let  $Q_{i,L}(w_L):=p_L(Q_i(w))$
be the $L$th projection image of $Q_i(w)$.
Then,  (\ref{eq:qsys5}) is equivalent to
a  series of equations ($L=1$, 2, \dots),
\begin{align}\label{eq:qsys6}
Q_{i,L}(w_L)+p_L(w_i) \prod_{j\in H}
(Q_{j,L}(w_L))^{G_{ij}}=1,
\end{align}
which  are further equivalent to
\begin{alignat}{2}
\label{eq:qsys16}
Q_{i,L}(w_L)&=1  & \quad &i\notin H_L,\\
\label{eq:qsys7}
Q_{i,L}(w_L)+w_i \prod_{j\in H_L}(Q_{j,L}(w_L))^{G_{ij}}&=1
& &i\in H_L.
\end{alignat}
Namely, a standard infinite $Q$-system is
 an infinite series of  standard finite $Q$-systems
which is compatible with the projections (\ref{eq:pl1}).
By  Proposition \ref{prop:qsys1},
(\ref{eq:qsys7}) uniquely determines  
$Q_{i,L}(w_L)$ for $i\in H_L$.
Furthermore, so determined
$(Q_{i,L}(w_L))_{L=1}^\infty$
belongs to $\bC[[w]]$, again because of the uniqueness
of the solution of (\ref{eq:qsys7}).
Therefore,

\begin{prop}\label{prop:iq1}
There exists a unique solution
$(Q_{i}(w))_{i\in H}$ of the standard
$Q$-system (\ref{eq:qsys5}),
whose $L$th projections  $Q_{i,L}(w_L):=p_L(Q_i(w))$
are determined by (\ref{eq:qsys16}) and (\ref{eq:qsys7}).
\end{prop}

\subsubsection{Canonical solution}

As we have seen in Example \ref{exmp:inf1},
the uniqueness property does not hold for a
general infinite $Q$-system (\ref{eq:qsys8}).
This is because, unlike the standard case,
the $L$th projection of (\ref{eq:qsys8})
is not necessarily a finite $Q$-system.
The non-uniqueness property 
also implies that, unlike the finite case,
(\ref{eq:qsys8})  does not
always reduce to
the standard one
\begin{align}\label{eq:qsys15}
Q'_i(w)+w_i \prod_{j\in H}(Q'_j(w))^{G'_{ij}}=1,
\quad
G'=GD^{-1}.
\end{align}
In fact,
 the relations (\ref{eq:qd1})
and (\ref{eq:qd8}) are no longer equivalent  due to
the infinite products therein.
However, the construction 
of a solution of a general $Q$-system
{}from a standard one in Theorem
\ref{thm:qunique1} still works.
We call the so obtained solution as the  {\em canonical
solution}.
Let us give a more intrinsic definition, however.

\begin{defn}
\label{def:can1}
We say that a solution $(Q_{i}(w))_{i\in H}$ of 
the $Q$-system (\ref{eq:qsys8}) is
{\em  canonical\/} if
it satisfies the following condition:
\par
(Inversion property): For any $i \in H$,
\begin{align} 
\label{eq:q22}
\prod_{j\in H}
\Bigl\{
\prod_{k\in H}
(Q_{k}(w))^{(D^{-1})_{ij}D_{jk}}
\Bigr\}
&=
Q_{i}(w).
\end{align}
\end{defn}

\begin{rem}
The condition (\ref{eq:q22}) is not trivial,
because, in general, one cannot freely
exchange the order of the
infinite double product therein.
\end{rem}

\begin{thm}\label{thm:cs4}
There exists a unique 
canonical solution of
the $Q$-system (\ref{eq:qsys8}),
which is given by
\begin{align}\label{eq:qd7}
Q_{i}(w)&=\prod_{j\in H}
(Q'_{j}(w))^{(D^{-1})_{ij}},
\end{align}
where $(Q'_i(w))_{i\in H}$ is
the unique solution of the standard $Q$-system (\ref{eq:qsys15}).
\end{thm}

\begin{proof}
First, we remark that the infinite product
 (\ref{eq:qd7}) exists,
because
its $L$th projection image  reduces to
the finite product
\begin{align}\label{eq:qd3}
Q_{i,L}(w_L)=\prod_{j\in H_L}
(Q'_{j,L}(w_L))^{(D^{-1})_{ij}}
\end{align}
due to (\ref{eq:qsys16}).
Let us show that
the family $(Q_i(w))_{i\in H}$ 
in (\ref{eq:qd7}) is a solution of
the $Q$-system (\ref{eq:qsys8}).
With the substitution of (\ref{eq:qd3}),
 the $L$th projection image of the first term 
in the LHS of (\ref{eq:qsys8}) is 
\begin{align}
\label{eq:qq2}
\begin{split}
\prod_{j\in H}(Q_{j,L}(w_L))^{D_{ij}}
&=
\prod_{j\in H}\Bigl\{
\prod_{k\in H_L}
(Q'_{k,L}(w_L))^{D_{ij}(D^{-1})_{jk}}
\Bigr\}\\
&=
\begin{cases}
Q'_{i,L}(w_L) &  i\in H_L\\
1=Q'_{i,L}(w_L) &  i\notin  H_L.\\
\end{cases}
\end{split}
\end{align}
In the second equality above, we exchanged
 the order of the products.
It is allowed because the double product
is a finite one (cf.\ Remark \ref{rem:dg1}).
The second term in the LHS of (\ref{eq:qsys8})
can be calculated in a similar way as follows:
\begin{align}
\label{eq:qq3}
\begin{split}
 \prod_{j\in H}
(Q_{j,L}(w_L))^{G_{ij}}
&=
\prod_{j\in H}\Bigl\{
\prod_{k\in H_L}
(Q'_{k,L}(w_L))^{G_{ij}(D^{-1})_{jk}}
\Bigr\}\\
&=
\prod_{k\in H_L}
(Q'_{k,L}(w_L))^{G'_{ik}}.\\
\end{split}
\end{align}
{}From (\ref{eq:qq2}) and (\ref{eq:qq3}),
we conclude that (\ref{eq:qsys8}) reduces to (\ref{eq:qsys15}).
Furthermore, by (\ref{eq:qq2}), we have
\begin{align}
\label{eq:qq8}
\prod_{j\in H}(Q_j(w))^{D_{ij}}
=
Q'_i(w).
\end{align}
Then, substituting (\ref{eq:qq8}) in
(\ref{eq:qd7}),  we obtain (\ref{eq:q22}).
Therefore, $(Q_i(w))_{i\in H}$ is a canonical solution
of (\ref{eq:qsys8}).
Next, we show the uniqueness.
Suppose that 
$(Q_{i}(w))_{i\in H}$ is a canonical solution of (\ref{eq:qsys8}).
We define $Q'_{i}(w)$ as
\begin{align}
\label{eq:q23}
Q'_{i}(w)=\prod_{j\in H}
(Q_{j}(w))^{D_{ij}}.
\end{align}
Then, by the inversion property (\ref{eq:q22}),
we have
\begin{align}\label{eq:qqd1}
Q_{i}(w)&=\prod_{j\in H}
(Q'_{j}(w))^{(D^{-1})_{ij}}.
\end{align}
Also, by (\ref{eq:qsys8}),
\begin{align}
\label{eq:qd10}
Q'_{i,L}(w_L)=1,
\quad
i\notin  H_L.
\end{align}
With (\ref{eq:qqd1}) and (\ref{eq:qd10}),
the same calculation as (\ref{eq:qq3})
shows that $(Q'_i(w))_{i\in H}$ is the (unique)
solution of
(\ref{eq:qsys15}).
Therefore, by (\ref{eq:qqd1}), 
$Q_i(w)$ is unique.
\end{proof}

\begin{exmp}
\label{exmp:inf2}
Let us find the canonical solution
of the $Q$-system (\ref{eq:qi2}) in Example \ref{exmp:inf1}.
We have
\begin{align}
\label{eq:dd1}
D_{ij}=-2\delta_{ij}+\delta_{i,j-1}+\delta_{i,j+1},
\quad
(D^{-1})_{ij}=-\min(i,j).
\end{align}
Let $H_L=\{1,\dots,L\}$.
By (\ref{eq:q23}) and (\ref{eq:qd10}),
the $L$th projection of the LHS of (\ref{eq:q22}) equals to
\begin{align}
\label{eq:qlim2}
\begin{split}
& \phantom{=}\ \prod_{j=1}^L \Bigl\{
\prod_{k=j-1}^{j+1}(Q_{k,L}(w_L))^{(D^{-1})_{ij}
D_{jk}}\Bigr\}\\
&=\Bigl(\prod_{k=1}^L \Bigl\{
\prod_{j=k-1}^{k+1}(Q_{k,L}(w_L))^{(D^{-1})_{ij}
D_{jk}}\Bigr\}\Bigr)\\
&\qquad 
\times (Q_{L+1,L}(w_L))^{(D^{-1})_{iL}D_{L,L+1}}
(Q_{L,L}(w_L))^{-(D^{-1})_{i,L+1}D_{L+1,L}}\\
&=
\Bigl(\prod_{k=1}^L (Q_{k,L}(w_L))^{\delta_{ik}}\Bigr)
(Q_{L+1,L}(w_L))^{-\min(i,L)}
(Q_{L,L}(w_L))^{\min(i,L+1)}\\
\end{split}
\end{align}
Therefore, the condition (\ref{eq:q22}) reads as
\begin{align}
Q_{i,L}(w_L)
&
=
\begin{cases}
Q_{i,L}(w_L)(Q_{L,L}(w_L)/Q_{L+1,L}(w_L))^i
& i\leq L\\
Q_{L,L}(w_L)(Q_{L,L}(w_L)/Q_{L+1,L}(w_L))^L
& i \geq L+1.\\
\end{cases}
\end{align}
This is equivalent to
\begin{align}
\label{eq:qi3}
Q_{i,L}(w_L)=Q_{L,L}(w_L),
\quad i\geq L+1.
\end{align}
Using (\ref{eq:qi1}) and (\ref{eq:qi3}),
one can easily obtain
\begin{align}
Q_1(w)=\prod_{j=1}^\infty (1-w_j)^{-1}.
\end{align}
Therefore, 
the canonical solution of (\ref{eq:qi2}) is given by
\begin{align}
Q_i(w)=\prod_{j=1}^\infty
(1-w_j)^{-\min(i,j)}.
\end{align}
\end{exmp}

\subsection{Power series formula}
\label{subsec:psf2}

Let $(Q_i(w))_{i\in H}$ be the
canonical solution of (\ref{eq:qsys8}),
and $(Q'_i(w))_{i\in H}$ be the
unique solution of the  standard
$Q$-system (\ref{eq:qsys15}).
For the matrix $D$ in (\ref{eq:qsys8}), 
let
$\bsnu(D)$
be the set of all
$\nu=(\nu_i)_{i\in H}$
such that $\nu_i\in \bC$ and,
for each $i$, the sum
$\sum_{j\in H} \nu_j (D^{-1})_{ji}$
exists
({\em i.e.,} 
 all but finitely many
 $\nu_j (D^{-1})_{ji}$ ($j\in H$) are zero).
For each $\nu\in \bsnu(D)$, we define
\begin{align}
\label{eq:qpow1}
Q^\nu_{D,G}(w)
:=\prod_{i\in H}
(Q_i(w))^{\nu_i}
=\prod_{i\in H}
\Bigl\{
\prod_{j\in H}
(Q'_{j}(w))^{\nu_i(D^{-1})_{ij}}
\Bigr\}.
\end{align} 
The last infinite product
 exists, because
its $L$th projection image
reduces to a finite product
 due to (\ref{eq:qsys16}) and
the definition of $\bsnu(D)$.
For each $\nu\in \bsnu(D)$,
let $\nu'=(\nu'_i)\in \bsnu(I)$,
$\nu'_i=\sum_{j\in H}\nu_j (D^{-1})_{ji}$.
Then, by (\ref{eq:qpow1}), we have
\begin{align}
\label{eq:qq11}
Q^{\nu}_{D,G}(w)&=
Q^{\nu'}_{I,G'}(w),
\quad G'=GD^{-1}.
\end{align}
It follows from (\ref{eq:qsys16})
and (\ref{eq:qq11}) that

\begin{lem}
\label{lem:pq1}
\begin{align}
p_L(Q^\nu_{D,G}(w))=
Q^{\nu'_L}_{I_L,G'_L}(w_L),
\end{align}
where the RHS is for the solution of the
finite $Q$-system
with the finite index set $H_L$,
and  
$I_L=(\delta_{ij})_{i,j\in H_L}$,
$G'_L=(G'_{ij})_{i,j\in H_L}$,
$\nu'_L=(\nu'_i)_{i\in H_L}$
are the $H_L$-truncations of $I$,  $G'$, $\nu'$,
respectively.
\end{lem}

For $D$, $G$ in (\ref{eq:qsys8}) and $\nu\in \bsnu(D)$,
we define the power series
$K^\nu_{D,G}(w)$ 
and
$R^\nu_{D,G}(w)$ 
by the superficially identical
formulae (\ref{eq:cKdef1})--(\ref{eq:fdef1})
with $D$, $G$, $\nu$, $N$, {\em etc.}, therein
being replaced by the ones for the infinite index set
$H$.

\begin{thm}[Power series formulae]
For the canonical solution $(Q_i(w))_{i\in H}$ of 
(\ref{eq:qsys8}) and  $\nu\in \bsnu(D)$,
let $Q^\nu_{D,G}(w)$ be the series in (\ref{eq:qpow1}).
Then,
\label{thm:qkr3}
\begin{align}
\label{eq:qkr2}
Q^\nu_{D,G}(w)=K^\nu_{D,G}(w)/K^0_{D,G}(w)=R^\nu_{D,G}(w).
\end{align}
\end{thm}

\begin{proof}
By Theorem \ref{thm:fqmain1} and Lemma \ref{lem:pq1},
it is enough to show that
\begin{align}
\label{eq:kl1}
p_L(K_{D,G}^\nu(w))=K_{I_L,G'_L}^{\nu'_L}(w_L),
\quad
p_L(R_{D,G}^\nu(w))=R_{I_L,G'_L}^{\nu'_L}(w_L).
\end{align}
By (\ref{eq:f1})--(\ref{eq:nldef1}),
(\ref{eq:kl1}) further reduces to
the following equality:
\begin{align}
\label{eq:pp1}
P_i(D,G;\nu,N)
=P_i(I_L,G'_L;\nu'_L,N_L),
\quad
N\in \cN_L,\
i\in H_L,
\end{align}
where $N_L=(N_i)_{i\in H_L}$ is the $H_L$-truncation
of $N$. 
\end{proof}

\section{$Q$-systems of KR type}
\label{sec:spec}

In this section,
we introduce a  class of infinite
$Q$-systems
which we call the $Q$-systems of KR type.
This is a  preliminary step towards 
the reformulation  of 
Conjecture \ref{conj:KR1}.

\subsection{Specialized  $Q$-systems}
\label{subsec:spec1}

Throughout the section,
 we take  the countable
index set as
\begin{align}
\label{eq:hinf1}
H=\{1,\dots,n\}\times \bN
\end{align}
for a given natural number $n$.
We choose the increasing sequence
$H_1\subset H_2\subset \cdots \subset H$
with $\varinjlim H_L=H$
as
$H_L=
\{1,\dots,n\}\times \{1,\dots,L\}$.
Let $y=(y_a)_{a=1}^n$ be a multivariable
 with $n$ components.

\begin{defn}
The following system of equations for
a family $(\cQ\am(y))_{(a,m)\in H}$ of
 power series of $y$ with the unit constant terms
is called a  {\em specialized (infinite)
$Q$-system\/}: For each $(a,m)\in H$,
\begin{align}\label{eq:qsys11}
\prod_{(b,k)\in H}
(\cQ\bk(y))^{D_{am,bk}}
+(y_a)^m \prod_{(b,k)\in H}(\cQ\bk(y))^{G_{am,bk}}=1,
\end{align}
where the infinite-size complex matrices
$D=(D_{am,bk})_{(a,m),(b,k)\in H}$
 and $G=(G_{am,bk})_{(a,m),(b,k)\in H}$
satisfy the same conditions (D) and (G') as in 
Definition \ref{def:qsys1}. 
A solution of (\ref{eq:qsys11}) is called
{\em canonical\/} if it satisfies
the condition
\begin{align} 
\label{eq:q24}
\prod_{(b,k)\in H}
\Bigl\{
\prod_{(c,j)\in H}
(\cQ\cj(y))^{(D^{-1})_{am,bk}D_{bk,cj}}
\Bigr\}
=
\cQ\am(y).
\end{align}
\end{defn}

Let $\bC[[y]]$ be the field of
power series of $y$ with the standard topology,
$J_L$ be the ideal of $\bC[[y]]$
generated by $(y_a)^{L+1}$'s ($a=1,\ \dots,\ n$),
and $\bC[[y]]_L$ be the quotient $\bC[[y]]/J_L$.
We can identify $\bC[[y]]$ with the projective limit of
the projective system,
\begin{align}\label{eq:proj2}
\bC[[y]]_1
\leftarrow \bC[[y]]_2\leftarrow \bC[[y]]_3
\leftarrow\cdots.
\end{align}
Let $w=(w\am)_{(a,m)\in H}$ be a multivariable,
and let $w(y)$ be the map with
\begin{align}
\label{eq:wy3}
w\am(y)&=(y_a)^m.
\end{align}
The map (\ref{eq:wy3}) induces the maps
$\psi_L$ and $\psi$ such that
\begin{align}
\label{eq:wycom1}
\begin{matrix}
 \bC[[w_L]]  
&  \leftarrow &\bC[[w]]\\
{\psi}_L \downarrow \ \ \ \  
 & &{\psi}\downarrow\  \\
\bC[[y]]_L &
 \leftarrow &\ \bC[[y]]. &\\
\end{matrix}
\end{align}
We call the image $\psi(f(w))\in \bC[[y]]$ 
the {\em specialization of $f(w)$},
and write it as  $f(w(y))$.
Explicitly, for $f(w)$ in (\ref{eq:f1}),
\begin{align}
f(w(y)) 
&= 
\sum_{M_1,\dots,M_n=0}^\infty
\Bigl(
\sum_{\genfrac{}{}{0pt}{1}{N\in \cN}
{ \sum_{m=1}^\infty m N\am=M_a}}
a_N
\Bigr)
\prod_{a=1}^n (y_a)^{M_a}.
\end{align}

\begin{thm}
\label{thm:cs2}
There exists a unique
canonical solution
of the specialized $Q$-system (\ref{eq:qsys11}),
which  is given by the specialization
$\cQ\am(y)=Q\am(w(y))$
of the canonical solution  $(Q\am(w))_{(a,m)\in H}$ 
 of the following $Q$-system:
\begin{align}\label{eq:qsys10}
\prod_{(b,k)\in H}
(Q\bk(w))^{D_{am,bk}}+w\am \prod_{(b,k)\in H}
(Q\bk(w))^{G_{am,bk}}=1.
\end{align}
\end{thm}
\begin{proof}
Since the map $\psi$ is
continuous,
it preserves the infinite product.
Therefore, 
the specialization of the
canonical solution of (\ref{eq:qsys10})
gives a canonical
solution of (\ref{eq:qsys11}).
Let us show the uniqueness.
By repeating the same proof for
Theorem \ref{thm:cs4},
the uniqueness is reduced 
to the one
for the standard case $D=I$.
Let us write
(\ref{eq:qsys11}) for $D=I$ as ($L=1,2,\dots$)
\begin{alignat}{2}
\cQ\am(y)\equiv &1 \mod J_L
&\quad & (a,m)\notin H_L,\\
\label{eq:qsys24}
\cQ\am(y)+(y_a)^m \prod_{(b,k)\in H_L}
(\cQ\bk(y))^{G_{am,bk}}\equiv &1 \mod J_L
& & (a,m)\in H_L.
\end{alignat}
These equations uniquely
determine $\cQ\am(y)$ mod $J_L$.
Since $L$ is arbitrary, $\cQ\am(y)$ is unique.
\end{proof}

By the specialization of Theorem \ref{thm:qkr3},
we immediately obtain
\begin{thm}[Power series formulae]
\label{thm:psf5}
Let $(\cQ\am(y))_{(a,m)\in H}$ be the canonical solution
of the $Q$-system  (\ref{eq:qsys11}).
Let $\cQ^\nu_{D,G}(y)=\prod_{(a,m)\in H}
(\cQ\am(w))^{\nu\am}$,
$\nu\in \bsnu(D)$.
Then,
\begin{align}
\label{eq:qkr8}
\cQ^\nu_{D,G}(y)=\cK^\nu_{D,G}(y)/\cK^0_{D,G}(y)=\cR^\nu_{D,G}(y),
\end{align}
where the series $\cK^\nu_{D,G}(y)=K^\nu_{D,G}(w(y))$
and $\cR^\nu_{D,G}(y)=R^\nu_{D,G}(w(y))$
are the specializations of the series
 in Theorem \ref{thm:qkr3}.
\end{thm}

\subsection{Convergence property}

Let us consider the special case of the
specialized $Q$-system
(\ref{eq:qsys11})
where the matrix $D$ and its inverse $D^{-1}$ are given by 
\begin{align}
\label{eq:krd1}
D_{am,bk}&=-\delta_{ab}(2\delta_{mk}-\delta_{m,k+1}
-\delta_{m,k-1}),\\
\label{eq:krd2}
(D^{-1})_{am,bk}&= -\delta_{ab}\min(m,k).
\end{align}
Then, (\ref{eq:qsys11}) is written in the form
($\cQ\azero(y)=1$)
\begin{align}\label{eq:qsys14}
\begin{split}
(\cQ\am(y))^2
&=
\cQ\amm(y)\cQ\amp(y)\\
&\qquad\qquad
+(y_a)^m (\cQ\am(y))^2\prod_{(b,k)\in H}
(\cQ\bk(y))^{G_{am,bk}}.
\end{split}
\end{align}

\begin{prop}
\label{prop:conv1}
A solution $(\cQ\am(y))$ of the specialized $Q$-system
 (\ref{eq:qsys14})
is canonical
if and only if it satisfies the following condition:
\par
(Convergence property):
\begin{align}
\label{eq:con1}
\text{For each $a$,
the limit\/ $\lim_{m\to\infty}\cQ\am(y)$ exists in
$\bC[[y]]$.} 
\end{align}
\end{prop}
\begin{proof}
Let $(\cQ\am(y))_{(a,m)\in H}$ be a solution of (\ref{eq:qsys14}).
The same calculation as (\ref{eq:qlim2}) in
Example \ref{exmp:inf2} shows that 
(\ref{eq:q24})
is equivalent to the following equality for each $L$
(cf.\  (\ref{eq:qi3})):
\begin{align}
\label{eq:qmod2}
\cQ\am(y)\equiv \cQ^{(a)}_L(y)
\mod J_L,
\quad m \geq L+1.
\end{align}
Clearly, the condition (\ref{eq:con1}) follows from 
the condition (\ref{eq:qmod2}).
Conversely,
assume the condition (\ref{eq:con1}).
By (\ref{eq:qsys14}), we have
\begin{align}
\label{eq:qmod1}
\cQ\am(y)/\cQ\amm(y)
\equiv
\cQ\amp(y)/\cQ\am(y)
\mod J_L,
\quad (m\geq L+1).
\end{align}
Because of (\ref{eq:con1}),
the both hand sides of (\ref{eq:qmod1})
are 1 mod $J_L$.
Thus,
 we have $\cQ\am(y)\equiv \cQ\amm(y)
\mod J_L$ ($m\geq L+1$).
Therefore,
(\ref{eq:qmod2}) holds.
\end{proof}

\subsection{$Q$-system of KR type
and denominator formula}

\begin{defn}
\label{defn:KR1}
A specialized 
$Q$-system (\ref{eq:qsys11})
is called
a {\em $Q$-system of
KR (Kirillov-Reshetikhin) type\/} if
the matrices $D$ and $G$ further satisfy the following
conditions:
\par
(KR-I) 
The matrix $D$ and its inverse $D^{-1}$ are given by 
(\ref{eq:krd1}) and (\ref{eq:krd2}).
\par
(KR-II)
There exists a well-order $\prec$ in $H$ 
such that
$G'=GD^{-1}$ has the form
\begin{align}
\label{eq:gd1}
G'_{am,bk}=g_{ab}m\quad \text{for $(a,m)\preceq (b,k)$},
\end{align}
where  $g_{ab}$ ($a,b=1,\dots,n$) are integers
with $\det_{1\leq a,b\leq n}
g_{ab}\neq 0$.
\end{defn}

\begin{exmp}
\label{exmp:g1}
Let $t_a>0$ 
and $h_{ab}$
($a,b=1,\dots,n$) be real numbers such that
$g_{ab}:=h_{ab}t_b$ are integers
 and $\det h_{ab}\neq 0$.
We define a well-order $\prec$ in $H$ as
follows: $(a,m)\prec (b,k)$
if $t_b m < t_a k$, or if
$t_b m = t_a k$ and $a<b$.
Then,
\begin{align}
\label{eq:gh1}
G'_{am,bk}=h_{ab}\min(t_b m,t_a k)
\end{align}
satisfies the condition (KR-II) with
$g_{ab}=h_{ab}t_b$.
\end{exmp}

Let $x=(x_a)_{a=1}^n$ be a multivariable
with $n$ components,
and $y(x)$ be the map
\begin{align}
\label{eq:yx2}
y_a(x)&=\prod_{b=1}^n (x_b)^{-g_{ab}},
\end{align}
where $g_{ab}$ are the integers in
(\ref{eq:gd1}).
We set
\begin{align}
\label{eq:bq1}
\bQ\am(x):=(x_a)^m \cQ\am(y(x)),
\end{align}
which are Laurent series of $x$.

\begin{prop}
The family $(\bQ\am(x))_{(a,m)\in H}$
satisfies a  system of 
equations $(\bQ\azero(x)=1)$,
\begin{align}\label{eq:qsys25}
\begin{split}
(\bQ\am(x))^2
&=
\bQ\amm(x)\bQ\amp(x)
+(\bQ\am(x))^2\prod_{(b,k)\in H}
(\bQ\bk(x))^{G_{am,bk}}.
\end{split}
\end{align}
\end{prop}
\begin{proof}
By comparing (\ref{eq:qsys14}) and (\ref{eq:qsys25}),
it is enough to prove the equality
\begin{align}
\label{eq:gsum1}
\sum_{k=1}^\infty G_{am,bk}(-k)=g_{ab}m.
\end{align}
{}Due to the condition (KR-II),
for given $(a,m)$ and $b$, there is some number $L$ such that
$G'_{am,bk}=g_{ab}m$ 
holds for any $k\geq L$.
Then, for $k>L$, we have
\begin{align}
G_{am,bk}&=\sum_{j=1}^\infty
G'_{am,bj}D_{bj,bk}=g_{ab}m(-2+1+1)=0.
\end{align}
Therefore, the LHS of (\ref{eq:gsum1}) is evaluated as
\begin{align}
\sum_{k=1}^{L} \sum_{j=1}^\infty
G'_{am,bj}D_{bj,bk}(-k)
=(L+1)G'_{am,bL}-L G'_{am,bL+1}=g_{ab} m.
\end{align}
\end{proof}

\begin{rem}
The relation (\ref{eq:qsys25}) is the original
 form of the $Q$-system  in
\cite{K2,K3,KR}, where the matrix $G$ is taken as
(\ref{eq:bformula2}).
See also (\ref{eq:gd2}) and (\ref{eq:bformula1}).
Note that, in the second term of the RHS in (\ref{eq:qsys25}),
the factor $(\cQ\am(y))^2$ is cancelled
by the factor in the product
for $(b,k)=(a,m)$, because  $G_{am,am}=-2$.
\end{rem}

\begin{prop}[Denominator formula]
\label{prop:denom11}
Let  $(\cQ\am(y))_{(a,m)\in H}$
be the canonical solution
of the $Q$-system of KR type (\ref{eq:qsys14}).
Let $\bK^0_{D,G}(x):=
\cK^0_{D,G}(y(x))$, where
$\cK^0_{D,G}(y)$ is
the power series  in (\ref{eq:qkr8}).
Then, the formula
\begin{align}
\label{eq:k07}
\bK^0_{D,G}(x)=
\det_{1\leq a,b\leq n}
\Bigl(
\frac{\partial \bQ\sa_1}{\partial x_b}(x)
\Bigr)
\end{align}
holds.
\end{prop}
A proof of Proposition
\ref{prop:denom11} is  given in Appendix \ref{sec:denom}.
In Conjecture \ref{conj:KR4},
Proposition \ref{prop:denom11}
will be used to
 identify $\bK^0_{D,G}(x)$ for some $G$
with the {\em Weyl denominators}
of the simple Lie algebras.

\section{$Q$-systems  and the Kirillov-Reshetikhin conjecture}
\label{sec:KR}

In this section, we reformulate    Conjecture
\ref{conj:KR1} in terms of
the canonical solutions of
certain
$Q$-systems of KR type (Conjecture \ref{conj:KR2}).
Then, we present several character formulae,
{\em all of which are equivalent to Conjecture \ref{conj:KR2}.}

\unitlength=0.8pt
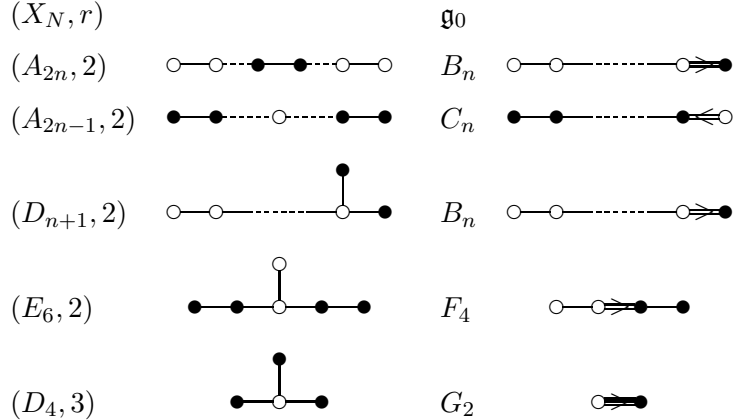
\begin{figure}[b]
\caption{The Dynkin diagrams of $X_N$ and $\fg_0$
for $r>1$.
The filled circles in $X_N$ correspond to
the ones in $\fg_0$ which are short roots of $\fg_0$.
}
\label{fig:d1}
\begin{tabular}[t]{lcclclrl}
$(X_N,r)$&&&$\fg_0$\\
$(A_{2n},2)$&
\begin{picture}(106,20)(-5,-5)
\multiput( 0,0)(20,0){2}{\circle{6}}
\multiput(40,0)(20,0){2}{\circle*{6}}
\multiput(80,0)(20,0){2}{\circle{6}}
\multiput( 3,0)(20,0){1}{\line(1,0){14}}
\multiput(23,0)(4,0){5}{\line(1,0){2}}
\multiput(43,0)(20,0){1}{\line(1,0){14}}
\multiput(59,0)(4,0){5}{\line(1,0){2}}
\multiput(83,0)(20,0){1}{\line(1,0){14}}
\end{picture}
&&$B_n$&
\begin{picture}(106,20)(-5,-5)
\multiput( 0,0)(20,0){2}{\circle{6}}
\multiput(80,0)(20,0){1}{\circle{6}}
\multiput(100,0)(20,0){1}{\circle*{6}}
\multiput( 3,0)(20,0){2}{\line(1,0){14}}
\multiput(63,0)(20,0){1}{\line(1,0){14}}
\multiput(82.85,-1)(0,2){2}{\line(1,0){14.3}}
\multiput(39,0)(4,0){6}{\line(1,0){2}}
\put(90,0){\makebox(0,0){$>$}}
\end{picture}
\\
$(A_{2n-1},2)$&
\begin{picture}(106,20)(-5,-5)
\multiput( 0,0)(20,0){2}{\circle*{6}}
\multiput(50,0)(20,0){1}{\circle{6}}
\multiput(80,0)(20,0){2}{\circle*{6}}
\multiput( 3,0)(20,0){1}{\line(1,0){14}}
\multiput(83,0)(20,0){1}{\line(1,0){14}}
\multiput(23,0)(4,0){6}{\line(1,0){2}}
\multiput(53,0)(4,0){6}{\line(1,0){2}}
\end{picture}
&&$C_n$&
\begin{picture}(106,20)(-5,-5)
\multiput( 0,0)(20,0){2}{\circle*{6}}
\multiput(80,0)(20,0){1}{\circle*{6}}
\multiput(100,0)(20,0){1}{\circle{6}}
\multiput( 3,0)(20,0){2}{\line(1,0){14}}
\multiput(63,0)(20,0){1}{\line(1,0){14}}
\multiput(82.85,-1)(0,2){2}{\line(1,0){14.3}}
\multiput(39,0)(4,0){6}{\line(1,0){2}}
\put(90,0){\makebox(0,0){$<$}}
\end{picture}
\\
$(D_{n+1},2)$&
\begin{picture}(106,40)(-5,-5)
\multiput( 0,0)(20,0){2}{\circle{6}}
\multiput(80,0)(20,0){1}{\circle{6}}
\multiput(100,0)(20,0){1}{\circle*{6}}
\put(80,20){\circle*{6}}
\multiput( 3,0)(20,0){2}{\line(1,0){14}}
\multiput(63,0)(20,0){2}{\line(1,0){14}}
\multiput(39,0)(4,0){6}{\line(1,0){2}}
\put(80,3){\line(0,1){14}}
\end{picture}
&&$B_n$&
\begin{picture}(106,20)(-5,-5)
\multiput( 0,0)(20,0){2}{\circle{6}}
\multiput(80,0)(20,0){1}{\circle{6}}
\multiput(100,0)(20,0){1}{\circle*{6}}
\multiput( 3,0)(20,0){2}{\line(1,0){14}}
\multiput(63,0)(20,0){1}{\line(1,0){14}}
\multiput(82.85,-1)(0,2){2}{\line(1,0){14.3}}
\multiput(39,0)(4,0){6}{\line(1,0){2}}
\put(90,0){\makebox(0,0){$>$}}
\end{picture}
\\
$(E_6,2)$&
\begin{picture}(86,40)(-5,-5)
\multiput(0,0)(20,0){2}{\circle*{6}}
\multiput(40,0)(20,0){1}{\circle{6}}
\multiput(60,0)(20,0){2}{\circle*{6}}
\put(40,20){\circle{6}}
\multiput(3,0)(20,0){4}{\line(1,0){14}}
\put(40, 3){\line(0,1){14}}
\end{picture}
&& $F_4$&
\begin{picture}(66,20)(-5,-5)
\multiput( 0,0)(20,0){2}{\circle{6}}
\multiput(40,0)(20,0){2}{\circle*{6}}
\multiput( 3,0)(40,0){2}{\line(1,0){14}}
\multiput(22.85,-1)(0,2){2}{\line(1,0){14.3}}
\put(30,0){\makebox(0,0){$>$}}
\end{picture}
\\
$(D_4,3)$&
\begin{picture}(46,40)(-5,-5)
\multiput(0,0)(20,0){1}{\circle*{6}}
\multiput(20,0)(20,0){1}{\circle{6}}
\multiput(40,0)(20,0){1}{\circle*{6}}
\put(20,20){\circle*{6}}
\multiput(3,0)(20,0){2}{\line(1,0){14}}
\put(20, 3){\line(0,1){14}}
\end{picture}
&& $G_2$&
\begin{picture}(26,20)(-5,-5)
\multiput( 0, 0)(20,0){1}{\circle{6}}
\multiput( 20, 0)(20,0){1}{\circle*{6}}
\multiput(2.68,-1.5)(0,3){2}{\line(1,0){14.68}}
\put( 3, 0){\line(1,0){14}}
\put(10, 0){\makebox(0,0){$>$}}
\end{picture}
\\
\end{tabular}
%
\end{figure}


\subsection{Quantum affine algebras}

We  formulate
   Conjecture
\ref{conj:KR1}
in the following setting:
Firstly, we translate the conjecture for
the KR modules of the {\em (untwisted) quantum affine
algebra\/} $U_q(X_n^{(1)})$,
based on the widely-believed correspondence
between the finite-dimensional
modules of $Y(X_n)$ and $U_q(X_n^{(1)})$ (for
the simply-laced case, see \cite{V}).
Secondly,
 we also include the {\em twisted\/} quantum affine algebra
case, following \cite{HKOTT}.

First, we introduce some notations.
Let $\fg=X_N$ be a 
complex simple Lie algebra of rank $N$.
We fix a Dynkin diagram automorphism $\sigma$ of $\fg$
with $r=\ord\, \sigma$.
Let $\fg_0$ be the $\sigma$-invariant
subalgebra of $\fg$;
namely,
\begin{align}
\label{eq:g0}
\begin{tabular}{c|cccccc}
$\fg$ &
$X_n$ & $A_{2n}$ & $A_{2n-1}$ & $D_{n+1}$ 
& $E_6$ & $D_4$ \\
$r$& 1& 2&2 & 2& 2&3\\
\hline
$\fg_0$ &
$X_n$ & $B_n$ & $C_n$ & $B_n$ & $F_4$ & $G_2$ 
\end{tabular}
\end{align}
See Figure \ref{fig:d1}.
Let $A'=(A'_{ij})$ ($i,j\in I$) and 
$A=(A_{ij})$ ($i,j\in I_\sigma$) be
the Cartan matrices of $\fg$ and $\fg_0$,
respectively,
where $I_\sigma$ is the set
of the $\sigma$-orbits on $I$.
We define the numbers $d'_i$, $d_i$, 
$\epsilon'_i$, $\epsilon_i$ ($i\in I$)
as follows:
$d'_i$ ($i\in I$) are coprime positive
integers such that $(d'_iA'_{ij})$ is symmetric;
$d_i$ ($i\in I_{\sigma}$) are coprime positive
integers such that $(d_iA_{ij})$ is symmetric,
and we set $d_i=d_{\pi(i)}$ ($i\in I$),
where $\pi:I\rightarrow I_\sigma$ is
the canonical projection;
$\epsilon'_i=r$ if $\sigma(i)=i$, and 1
otherwise;
$\epsilon_i=2$ if $A'_{i\sigma(i)}<0$, and 1
otherwise.
It immediately follows that
$d_i'=d_i$ and $\epsilon'_i=1$ if\ $r=1$;
$d'_i=1$ if $r>1$; $\epsilon_i=1$ if $X^{(r)}_N
\neq A^{(2)}_{2n}$. 
It is easy to check the following relations:
Set $\kappa_0=2$ if $X^{(r)}_N
= A^{(2)}_{2n}$, and 1 otherwise.
Then,
\begin{align}
\label{eq:car2}
\kappa_0 d'_i \sum_{s=1}^r
A'_{i\sigma^s(j)} &= d_i A_{\pi(i)\pi(j)},\\
\label{eq:car1}
\kappa_0 \epsilon'_i d'_i
&=\epsilon_i d_i.
\end{align}
For $q\in \bC^{\times}$,
we set  $q'_i=q^{\kappa_0 d'_i}$, $q_i=q^{d_i}$,
and $[n]_q=(q^n-q^{-n})/(q-q^{-1})$.

We use the ``second realization''
of the quantum affine algebra
$U_q=U_q(X^{(r)}_N)$ \cite{D2,J} with the
generators $X_{ik}^{\pm}$
($i \in I, k\in \bZ$),
$H_{ik}$
($i \in I, k\in \bZ\setminus \{0\}$),
$K_{i}^{\pm1}$
($i \in I$), and the central elements $c^{\pm1/2}$.
As far as  finite-dimensional $U_q$-modules
are concerned, we can set $c^{\pm1/2}=1$.
Some of the
defining relations
in the quotient (the {\em quantum loop
algebra\/}) $U_q/( c^{\pm1/2}-1)$
are presented
below to fix notations
(here we follow the convention in \cite{CP2,CP3}):

\begin{align}
\label{eq:rel1}
X^\pm_{\sigma(i)k}&=\omega^k X^\pm_{ik},
\quad
H_{\sigma(i)k}=\omega^k H_{ik},
\quad
K^{\pm1}_{\sigma(i)k}= K^{\pm1}_{ik},\\
\label{eq:rel2}
K_i X^\pm_{jk}K^{-1}_i
&=q^{\pm \kappa_0 d'_i \sum_{s=1}^r A'_{i\sigma^s(j)}}
X^\pm_{jk},\\
\label{eq:rel3}
[H_{ik},X^{\pm}_{jl}]
&
=\pm
\frac{1}{k}
\biggl(
\sum_{s=1}^r
[k \kappa_0 d'_i A'_{i\sigma^s(j)}]_{q}
\omega^{sk}\biggr)X^\pm_{j,k+l},
\\
\label{eq:rel4}
[X^+_{ik}, X^-_{jl}]&=
\biggl(
\sum_{s=1}^r
\delta_{\sigma^s(i),j}
\omega^{sl}
\biggr)
\frac{\Psi^+_{i,k+l}-\Psi^-_{i,k+l}}
{q_i-q^{-1}_i},
\end{align}
where $\omega=\exp(2\pi i/r)$, and
$\Psi^\pm_{ik}$ ($i\in I, k\in \bZ$)  are defined
by
\begin{align}
\label{eq:rel5}
\sum_{k=0}^\infty 
\Psi^\pm_{i,\pm k}u^{k}
=
K^{\pm1}_i
\exp
\biggl(
\pm(q_i-q^{-1}_i)
\sum_{l=1}^\infty
H_{i,\pm l}u^{l}
\biggr)
\end{align}
with $\Psi^\pm_{ik}=0$ ($\pm k<0$).

\begin{rem}
In \cite{CP3}, there are some misprints
which are relevant here.
Namely, the relation $[H_{ik},X^\pm_{jl}]$ should read
(\ref{eq:rel3}) here; in Proposition 2.2 and 
Theorem 3.1 (ii), $q$ should read $q_i$ for such
$i$ that
$\sigma(i)\neq i$ and $a_{i\sigma(i)}\neq 0$ therein.
We thank V.\ Chari for the correspondence concerning these
points.
\end{rem}

\par
Let $V(\psi^\pm)$ denote
the irreducible $U_q$-module
with a highest weight vector $v$
and the highest weight $\psi^\pm=
(\psi^\pm_{ik})$, namely,
\begin{align}
\label{eq:hw1}
X^+_{ik}v&=0,\\
\label{eq:hw2}
\Psi^\pm_{ik}v&=\psi^\pm_{ik}v,
\quad \psi^\pm_{ik}\in \bC.
\end{align}

The following theorem gives the classification
of the finite-dimensional $U_q$-modules:
\begin{thm}[Theorem 3.3 \cite{CP2}, Theorem 3.1 \cite{CP3}]
The $U_q(X^{(r)}_N)$-module $V(\psi^\pm)$ is
finite-dimensional if
and only if there exist $N$-tuple of polynomials
$(P_i(u))_{i\in I}$ with the unit constant terms
such that
\begin{align}
\label{eq:dp1}
\sum_{k=0}^\infty \psi^+_{ik}u^k
=
\sum_{k=0}^\infty \psi^-_{i,-k}u^{-k}
=
q'_i{}^{\epsilon'_i\deg P_i}
\frac{P_i(q'_i{}^{-2\epsilon'_i}
u^{\epsilon'_i})}{P_i(u^{\epsilon'_i})},
\end{align}
where the first two terms are the Laurent
expansions of the third term about
$u=0$ and $u=\infty$, respectively.
\end{thm} 

The  polynomials $(P_i(u))_{i\in I}$ are called
the {\em Drinfeld polynomials} of $V(\psi^\pm)$.
It follows from (\ref{eq:car1}),
(\ref{eq:hw2}), and (\ref{eq:dp1}) that
\begin{align}
\label{eq:dp2}
K_i^{\pm1} v =
q'_i{}^{\pm\epsilon'_i\deg P_i}v
=
q_i^{\pm\epsilon_i\deg P_i}v.
\end{align}

\subsection{The KR modules}

We take an inclusion $\iota:I_\sigma \hookrightarrow I$
such that $\pi\circ\iota=\mathrm{id}$,
and regard $I_\sigma$ as a subset of $I$.
Let us label the set
$I_\sigma$ with $\{1,\dots,n\}$.
The Drinfeld polynomials (\ref{eq:dp1}) satisfy 
the relation $P_{\sigma(i)}(u)=P_i(\omega u)$
($\sigma(i)\neq i$)
by (\ref{eq:rel1}) and (\ref{eq:rel5}).
Therefore, it is enough to specify the
polynomials $P_i(u)$ only for those $i\in 
\{1,\dots,n\}\subset I$.

We set $H=\{1,\dots,n\}\times \bN$ as in (\ref{eq:hinf1}).
\begin{defn}
For each $(a,m)\in H$ and $\zeta\in \bC^\times$,
let $W\am(\zeta)$ be the finite-dimensional
irreducible $U_q$-module whose  Drinfeld
polynomials $P_b(u)$ ($b=1,\dots,n$)
are specified as follows:
$P_b(u)=1$ for $b\neq a$,
and
\begin{align}
\label{eq:dp3}
P_a(u)=
\prod_{k=1}^m
(1-\zeta q'_a{}^{\epsilon'_a(m+2-2k)}u).
\end{align}
We call $W\am(\zeta)$ a
{\em KR (Kirillov-Reshetikhin) module}.
\end{defn}

By (\ref{eq:car2}) and (\ref{eq:rel2}),
we see that  $X^\pm_{a0}$ and $K^{\pm1}_a$
 ($a=1,\dots,n$)
generate the subalgebra $U_{q}(\fg_0)$.
It is well known that
all $W\am(\zeta)$ ($\zeta\in \bC^\times$) share the same
$U_{q}(\fg_0)$-module structure.
If we set $K^{\pm1}_a=q_a^{\pm H_a}$ and take the limit
$q\rightarrow 1$,
$X^\pm_{a0}$ and $H_a$ ($a=1,\dots,n$)
generate the Lie algebra $\fg_0$.
Accordingly,  $W\am(\zeta)$ is  equipped with
$\fg_0$-module structure.
We call its $\fg_0$-character the
 {\em $\fg_0$-character of $W\am(\zeta)$}.
The $\fg_0$-highest weight of $W\am(\zeta)$,
in the same sense as above, is
$m\epsilon_a\Lambda_a$
by (\ref{eq:dp2}) and (\ref{eq:dp3}).

\subsection{The Kirillov-Reshetikhin conjecture}
\label{subsec:KRc}

We define the matrix
$G'=(G'_{am,bk})_{(a,m),(b,k)\in H}$ with the entry
\begin{align}
\label{eq:gd2}
G'_{am,bk}&=
\sum_{s=1}^r\frac{d'_b}{\epsilon'_b}
 A'_{b\sigma^s(a)}
\min\Bigl(
\frac{m}{d'_b},\frac{k}{d'_a}\Bigr)
\\
\label{eq:gd3}
&=
\begin{cases}
d_b A_{ba}
\min(\frac{m}{d_b},
\frac{k}{d_a})& r=1\\
\frac{1}{\epsilon_b}A_{ba}
\min(m,k)& r>1.
\end{cases}
\end{align}
It follows from  (\ref{eq:gd3}) and Example \ref{exmp:g1}
that $G'$ satisfies the condition
(KR-II) in Definition
\ref{defn:KR1} with $g_{ab}=A_{ba}/\epsilon_b$.
Below, we consider the $Q$-system of KR type
with the matrix $G:=G'D$,
where $D$ is the matrix in (\ref{eq:krd1}).
By using (A.6) of \cite{KN2}),
the entry of $G$ is explicitly
written as
\begin{align}
\label{eq:bformula1}
G_{am,bk}=
\begin{cases}
-\frac{1}{\epsilon_b}
A_{ba}\delta_{m, k}
& r>1 \\
-A_{ba}
(\delta_{m,2k-1}+2\delta_{m,2k}+\delta_{m,2k+1})
&
d_b/d_a=2\\
-A_{ba}
(\delta_{m,3k-2}+2\delta_{m,3k-1}+3\delta_{m,3k}&
d_b/d_a=3\\
\qquad\qquad\qquad
+2\delta_{m,3k+1}+\delta_{m,3k+2})&\\
-A_{ab}\delta_{d_a m, d_b k}
& \text{otherwise}.
\end{cases}
\end{align}

Let $\alpha_a$ and $\Lambda_a$ ($a=1,\dots,n$)
be the simple roots and the fundamental weights
of $\fg_0$. We set
\begin{align}
\label{eq:ya0}
x_a&=e^{\epsilon_a\Lambda_a},
\quad y_a=e^{-\alpha_a}.
\end{align}
Then, they satisfy the relation (\ref{eq:yx2})
for the above $g_{ab}$;
namely,
\begin{align}
\label{eq:yx4}
y_a=\prod_{b=1}^n x_b^{-A_{ba}/\epsilon_b}.
\end{align}

\begin{defn}
\label{def:fc1}
Let $\bQ\am(x)$ be the Laurent polynomial of $x=(x_a)_{a=1}^n$
representing the $\fg_0$-character of the KR module
$W\am(\zeta)$.
Then, $\cQ\am(y):=(x_a)^{-m}\bQ\am(x)|_{x=x(y)}$,
where $x(y)$ is the inverse map of  (\ref{eq:yx4}),
is a polynomial of $y=(y_a)_{a=1}^n$
with the unit constant term.
We call $\cQ\am(y)$ the {\em normalized $\fg_0$-character}
of $W\am(\zeta)$.
\end{defn}

Now we present a reformulation of 
Conjecture \ref{conj:KR1}.
This is the main statement of the paper.

\begin{conj}
\label{conj:KR2}
Let $\cQ\am(y)$ be the normalized $\fg_0$-character
of the KR module $W\am(\zeta)$ of $U_q(X^{(r)}_N)$.
Then, the family $(\cQ\am(y))_{(a,m)\in H}$  is characterized
as the canonical
solution of the $Q$-system of KR type (\ref{eq:qsys14})
with $G$ given in (\ref{eq:bformula1}).
\end{conj}

Let $\cQ^{\nu}(y)
=\prod_{(a,m)} (\cQ\am(y))^{\nu\am}$
for $\nu\in \bsnu(D)$.
{}
By Theorem \ref{thm:psf5},
Conjecture \ref{conj:KR2} is equivalent to
\begin{conj}[\cite{KN2}]
\label{conj:KR3}
The formulae
\begin{align}
\label{eq:fc1}
\cQ^\nu(y)=\cK^\nu_{D,G}(y)/
\cK^0_{D,G}(y)=\cR^\nu_{D,G}(y)
\end{align}
hold, where $\cK^\nu_{D,G}(y)$
and $\cR^\nu_{D,G}(y)$ are the power series  in
(\ref{eq:qkr8}) with $D$ in (\ref{eq:krd1})
and $G$ in (\ref{eq:bformula1}).
Therefore, $\cR^\nu_{D,G}(y)$ is
a {\em polynomial\/} of $y$,
and its coefficients are identified with the $\fg_0$-weight
multiplicities of the tensor product
$\bigotimes_{(a,m)\in H}W\am(\zeta\am)^{\otimes \nu\am}$,
where $\zeta\am$ are arbitrary.
\end{conj}

\subsection{Equivalence to Conjecture \ref{conj:KR1}}
\label{subsec:Rc}

Let  $\Delta_+^{\fg}$ denote the set of all the positive
roots of $\fg$.
Originally,
 Conjecture \ref{conj:KR2} is formulated
for $X^{(r)}_N \neq A^{(2)}_{2n}$
as follows (cf.\ Conjecture \ref{conj:KR1}):

\begin{conj}[\cite{K1,K2,HKOTY,HKOTT}]
\label{conj:KR4}
For $X^{(r)}_N \neq A^{(2)}_{2n}$,
the formula
\begin{align}
\label{eq:krc3}
\cQ^\nu(y)=
\frac{
\cK^\nu_{D,G}(y)
}{\displaystyle
\prod_{\alpha \in \Delta_+^{\fg_0}}(1-e^{-\alpha})}
\end{align}
holds, where
$\cK^\nu_{D,G}(y)$ is the power series  in
(\ref{eq:qkr8}) with 
$D$ in (\ref{eq:krd1})
and $G$ in (\ref{eq:bformula1}).
Therefore, $\cK^\nu_{D,G}(y)$ is
a {\em polynomial\/} of $y$,
and its coefficients are identified with 
the multiplicities of the $\fg_0$-irreducible
components of the tensor product
$\bigotimes_{(a,m)\in H}W\am(\zeta\am)^{\otimes \nu\am}$,
where $\zeta\am$ are arbitrary.
\end{conj}

\par\noindent
{\em Proof of the equivalence between
Conjectures \ref{conj:KR3} and \ref{conj:KR4}
for $X^{(r)}_N \neq A^{(2)}_{2n}$.}
Suppose that  Conjecture \ref{conj:KR4} holds.
Then, setting $\nu=0$ in (\ref{eq:krc3}),
we have
\begin{align}
\label{eq:denom5}
\cK^0_{D,G}(y)=\prod_{\alpha\in \Delta_+^{\fg_0}}(1-e^{-\alpha}).
\end{align}
Therefore, $\cQ^\nu(y)=\cK^\nu_{D,G}(y)/
\cK^0_{D,G}(y)$ holds.
Conversely, suppose that 
the family of the normalized $\fg_0$-characters
$(\cQ^{(a)}_m(y))_{(a,m)\in H}$
is the
canonical solution of (\ref{eq:qsys14}).
Then, the equality (\ref{eq:denom5})
follows from Proposition \ref{prop:denom11}
and the lemma below.
\qed

\begin{lem}
\label{lem:denom1}
Let $\fg$ be a complex simple Lie algebra of rank $n$,
and $\alpha_a$ and $\Lambda_a$ be the
simple roots and the fundamental weights of
$\fg$.
We set $x_a=e^{\Lambda_a}$, $y_a=e^{-\alpha_a/k_a}$,
where $k_a$ ($a=1,\dots,n$) are 1 or 2.
Suppose that
$f_a(y)$ $(a=1,\dots,n)$ are polynomials of $y$
with the unit constant terms
such that  $\hf_a(x)=x_a f_a(y(x))$
are invariant under the action of the 
Weyl group of $\fg$.
Then,
\begin{align}
\det_{1\leq a,b\leq n}
\Bigl(
\frac{\partial \hf_a}{\partial x_b}(x)\Bigr)
=\prod_{\alpha\in \Delta_+^{\fg}}(1-e^{-\alpha}).
\end{align}
\end{lem}
\begin{proof}
The proof is the same as the one for Lemma 8.6 in
 \cite{HKOTY}.
\end{proof}

In the case $A^{(2)}_{2n}$,
(\ref{eq:denom5}) does not hold
under Conjecture \ref{conj:KR3},
because the assumption in
Lemma \ref{lem:denom1} is not satisfied
by (\ref{eq:ya0}).
We treat the case $A^{(2)}_{2n}$
separately below.

\subsection{The $A^{(2)}_{2n}$ case}
\label{subsec:a22}

\subsubsection{The $B_n$-character}

\unitlength=0.8pt
\begin{figure}[t]
\caption{The Dynkin diagram of $A^{(2)}_{2n}$.
The upper and lower labels respect
the subalgebra $B_n$ and $C_n$, respectively.
}
\label{fig:d2}
\begin{picture}(106,30)(-5,-7)
\multiput(2.85,-1)(0,2){2}{\line(1,0){14.3}}
\multiput(0,0)(20,0){2}{\circle{6}}
\multiput(80,0)(20,0){1}{\circle{6}}
\multiput(100,0)(20,0){1}{\circle{6}}
\multiput(23,0)(20,0){1}{\line(1,0){14}}
\multiput(63,0)(20,0){1}{\line(1,0){14}}
\multiput(82.85,-1)(0,2){2}{\line(1,0){14.3}}
\multiput(39,0)(4,0){6}{\line(1,0){2}}
\put(10,0){\makebox(0,0){$>$}}
\put(90,0){\makebox(0,0){$>$}}
\put(0,17){\makebox(0,0)[t]{$0$}}
\put(20,17){\makebox(0,0)[t]{$1$}}
\put(80,17){\makebox(0,0)[t]{$n\!\! -\!\! 1$}}
\put(100,14){\makebox(0,0)[t]{$n$}}
\put(0,-10){\makebox(0,0)[t]{$n$}}
\put(20,-7){\makebox(0,0)[t]{$n\!\! -\!\! 1$}}
\put(80,-7){\makebox(0,0)[t]{$1$}}
\put(100,-7){\makebox(0,0)[t]{$0$}}
\end{picture}
\end{figure}
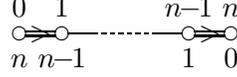

For $A^{(2)}_{2n}$, $\fg_0=B_n$.
Let $\{1,\dots,n\}$ label $I_\sigma$
as the upper label    in Figure \ref{fig:d2}.
Accordingly, $\epsilon_a=1$ for $a=1,\dots,n-1$,
and 2 for $a=n$.
We continue to set
$y_a=e^{-\alpha_a}$
as in Section \ref{subsec:KRc}.
We will show later, in (\ref{eq:denom9}) and 
(\ref{eq:denom7}), that
under Conjecture \ref{conj:KR2} the following formula holds 
instead of the formula (\ref{eq:denom5}):
\begin{align}
\label{eq:denom6}
\cK^0_{D,G}(y)=
\prod_{a=1}^n\Bigl(1+\prod_{k=a}^n y_k\Bigr)
\prod_{\alpha\in \Delta_+^{B_n}}(1-e^{-\alpha}).
\end{align}
Therefore,
Conjecture \ref{conj:KR2} 
for  $X^{(r)}_N = A^{(2)}_{2n}$ is equivalent to

\begin{conj}
\label{conj:KR5}
For  $X^{(r)}_N = A^{(2)}_{2n}$,
the formula 
\begin{align}
\label{eq:krc4}
\cQ^\nu(y)=
\frac{
\cK^\nu_{D,G}(y)
\prod_{a=1}^n\bigl(1+\prod_{k=a}^n y_k\bigr)^{-1}
}{\displaystyle
\prod_{\alpha \in \Delta_+^{B_n}}(1-e^{-\alpha})}
\end{align}
holds
for the normalized $B_n$-characters of the KR-modules.
\end{conj}

\subsubsection{The $C_n$-character}
\label{subsubsec:Cn}

As is well-known,
$U_q(A^{(2)}_{2n})$ has
a realization with the
``Chevalley generators''
$X^\pm_{a}$ and $K^{\pm1}_a$ ($a=0,\dots,n$)
 (e.g.\ \cite[Proposition 1.1]{CP3}).
Among them,
$X^\pm_{a}$ and $K^{\pm1}_a$ ($a=1,\dots,n$)
are identified with $X^\pm_{a0}$, $K^{\pm1}_a$
in (\ref{eq:rel1})--(\ref{eq:rel5}),
and generate the subalgebra
$U_{q}(B_n)$.
On the other hand,
$X^\pm_{a}$ and $K^{\pm1}_a$ ($a=0,\dots,n-1$)
generate the  subalgebra $U_{q^2}(C_n)$.
See Figure \ref{fig:d2}.
If we set $K_a=q_a^{H_a}$ 
($a=0,\dots,n-1$), where $q_0=q^{d_0}$, $d_0=4$,
then $X^\pm_{a}$ and $H_a$ ($a=0,\dots,n-1$)
generate the Lie algebra $C_n$
in the limit $q\rightarrow 1$.
This provides $W\am(\zeta)$ with
the $C_n$-module structure,
by which 
the {\em $C_n$-character of $W\am(\zeta)$}
is defined.

Let $\dot{\alpha}_a$ and $\dot{\Lambda}_a$
($a=1,\dots,n$) be the simple roots and
the fundamental weights labeled with
the lower label in Figure \ref{fig:d2}.
By looking at the same $U_q$-module
as $B_n$ and $C_n$-modules as above,
a linear bijection $\phi:\fh^*\rightarrow \dot{\fh}^{*}$
is induced,
where $\fh^*$ and $\dot{\fh}^{*}$ are the
duals of the Cartan subalgebras of
$B_n$ and $C_n$, respectively.

\begin{lem}
\label{lem:w1}
Under the bijection $\phi$,
we have the correspondence $(\dot{\Lambda}_0=0)$:
\begin{align}
\epsilon_a\Lambda_a &\mapsto \dot{\Lambda}_{n-a}-
\dot{\Lambda}_n,\\
\alpha_a &\mapsto 
\begin{cases}
\dot{\alpha}_{n-a}&
a=1,\dots,n-1\\
 -(\dot{\alpha}_1+\dots+
\dot{\alpha}_{n-1}+\frac{1}{2}\dot{\alpha}_n)
& a=n.\\
\end{cases}
\end{align}
\end{lem}
\begin{proof}
It is  obtained from the
relations among $H_i$ and $\alpha_i$ for $A^{(2)}_{2n}$
\cite{Kac}:
\begin{align}
0=c=\sum_{i=0}^n a^\vee_i H_i,
\quad
0=\delta=\sum_{i=0}^n a_i \alpha_i,
\end{align}
where $(a^\vee_0,\dots,a^\vee_n)=
(2,\dots,2,1)$ and 
$(a_0,\dots,a_n)=(1,2,\dots,2)$
for the upper label in Figure \ref{fig:d2}.
\end{proof}

Let $\cW(X_n)$ denote the Weyl group of $X_n$.

\begin{lem}
\label{lem:w2}
There is an element $s\in \cW(C_n)$ which acts
on $\dot{\fh}^*$ as follows:
\begin{align}
\phi(\epsilon_a\Lambda_a)
 &\mapsto \dot{\Lambda}_a
\quad (a=1,\dots,n),
\\
\phi(\alpha_{a}) &\mapsto 
\frac{1}{\epsilon_a}
\dot{\alpha}_a
\quad (a=1,\dots,n).
\end{align}
\end{lem}

\begin{proof}
We take the standard orthonormal basis $\varepsilon_a$
of $\dot{\fh}^*$. Let $s$ be the element
such that $s:\varepsilon_a\mapsto -\varepsilon_{n-a+1}$.
Then,
\begin{align}
\label{eq:lt1}
&\dot{\Lambda}_{n-a}-
\dot{\Lambda}_n
=-(\varepsilon_{n-a+1}+\dots+\varepsilon_n)
 \mapsto 
\varepsilon_{1}+\dots+\varepsilon_a=\dot{\Lambda}_a,\\
\label{eq:lt2}
&\dot{\alpha}_{n-a}
=\varepsilon_{n-a}-\varepsilon_{n-a+1}
 \mapsto 
\varepsilon_{a}-\varepsilon_{a+1}
=
\dot{\alpha}_a
 \quad (a=1,\dots,n-1),\\
\label{eq:lt3}
& -(\dot{\alpha}_1+\dots+
\dot{\alpha}_{n-1}+\frac{1}{2}\dot{\alpha}_n)
=-\varepsilon_1
\mapsto 
\varepsilon_n=
\frac{1}{2}\dot{\alpha}_n.
\end{align}
\end{proof}

According to (\ref{eq:lt1})--(\ref{eq:lt3}),
we set
\begin{align}
\label{eq:ya3}
x_a=e^{\dot{\Lambda}_a},\quad
y_a=e^{-\dot{\alpha}_a/\epsilon_a}.
\end{align}
Then, the relation (\ref{eq:yx4}) is preserved,
since  $\phi$ and $s$ above are linear.
Lemma  \ref{lem:w2}
assures that the following definition is well-defined.

\begin{defn}
\label{def:fc2}
Let ${\bQ}\am(x)$ be the Laurent polynomial of $x=(x_a)_{a=1}^n$
representing the $C_n$-character of the KR module
$W\am(\zeta)$.
Then, ${\cQ}\am(y):=(x_a)^{-m}{\bQ}\am(x)|_{x=x(y)}$
is a polynomial of $y=(y_a)_{a=1}^n$
with the unit constant term.
We call ${\cQ}\am(y)$ the {\em normalized $C_n$-character}
of  $W\am(\zeta)$.
\end{defn}
Moreover,
by Lemma  \ref{lem:w2} 
and the $\cW(C_n)$-invariance of the $C_n$-character
of $W\am(\zeta)$, we have
\begin{prop}
\label{prop:KR1}
 The normalized $B_n$-character and the
normalized $C_n$-character
of $W\am(\zeta)$ of $U_q(A^{(2)}_{2n})$ coincide
as polynomials of $y$.
\end{prop}

Thus, Conjecture \ref{conj:KR2}
for the normalized $B_n$-characters
of $A^{(2)}_{2n}$ is applied
for the normalized $C_n$-characters as well.
Furthermore, 
in contrast to the $B_n$ case,
Lemma \ref{lem:denom1} is now applicable for (\ref{eq:ya3}).
Therefore,
under Conjecture \ref{conj:KR2}, we have
\begin{align}
\label{eq:denom9}
\cK^0_{D,G}(y)=
\prod_{{\alpha} \in {\Delta}^{C_n}_+}(1-e^{-{\alpha}}).
\end{align}
Hence, we conclude that 
Conjecture \ref{conj:KR2} 
for  $X^{(r)}_N = A^{(2)}_{2n}$ is also equivalent to

\begin{conj}[\cite{HKOTT}]
\label{conj:KR6}
For  $X^{(r)}_N = A^{(2)}_{2n}$,
the formula
\begin{align}
\label{eq:krc5}
\cQ^\nu(y)=
\frac{
\cK^\nu_{D,G}(y)}{\displaystyle
\prod_{{\alpha} \in {\Delta}^{C_n}_+}(1-e^{-{\alpha}})}
\end{align}
holds
for the normalized $C_n$-characters of the KR-modules,
where $y$ is specified as (\ref{eq:ya3}).
\end{conj}

The following relation is easily derived from the
explicit expressions of the Weyl denominators
of $B_n$ and $C_n$ (e.g.\ \cite{FH}):
\begin{align}
\label{eq:denom7}
\prod_{\alpha\in {\Delta}^{C_n}_+}(1-e^{-\alpha})
=\prod_{a=1}^n\Bigl(1+\prod_{k=a}^n y_k\Bigr)
\prod_{\alpha\in \Delta^{B_n}_+}(1-e^{-\alpha}),
\end{align}
where the equality holds under the following identifications:
$y_a=e^{-\dot{\alpha}_a/\epsilon_a}$ for the LHS
and  $y_a=e^{-\alpha_a}$ for the RHS
under the label in Figure \ref{fig:d2}.
{}From (\ref{eq:denom9}) and (\ref{eq:denom7}),
we obtain (\ref{eq:denom6}).

\subsection{Characters for the rank $n$ subalgebras}

The procedure to deduce the $C_n$-characters
from the $B_n$-characters for $A^{(2)}_{2n}$
in Section \ref{subsec:a22}
is also applicable to the $\dot{\fg}$-characters
for any rank $n$ subalgebra $\dot{\fg}\neq \fg_0$ of $X^{(r)}_N$.
(The characters of the lower rank subalgebras are obtained 
by their specializations.)
Let us demonstrate how it works in two examples:
\par Case I. $X^{(r)}_N=B^{(1)}_n$, $\fg_0=B_n$, $\dot{\fg}=D_n$.
\par Case II. $X^{(r)}_N=A^{(2)}_{2n-1}$, $\fg_0=C_n$, $\dot{\fg}=D_n$.
\par
Let ${\alpha}_a$ and ${\Lambda}_a$
(resp.\ $\dot{\alpha}_a$ and $\dot{\Lambda}_a$)
($a=1,\dots,n$) be the simple roots and
the fundamental weights of $\fg_0$
(resp.\ $\dot{\fg}$)
labeled with
the upper (resp.\ lower) label in Figure \ref{fig:d3}.
As in Section \ref{subsec:a22},
a linear bijection $\phi:\fh^*\rightarrow \dot{\fh}^{*}$
is induced,
where $\fh^*$ and $\dot{\fh}^{*}$ are the
duals of the Cartan subalgebras of
$\fg_0$ and $\dot{\fg}$, respectively.

\unitlength=0.8pt
\begin{figure}[t]
\caption{The Dynkin diagrams of $B^{(1)}_n$
and $A^{(2)}_{2n-1}$.
The upper and lower labels respect
the subalgebra $\fg_0$ and $\dot{\fg}$, respectively.
}
\label{fig:d3}
\begin{picture}(106,85)(-5,-43)
\multiput(2.5,17.5)(0,2){1}{\line(1,-1){15}}
\multiput(2.5,-17.5)(0,2){1}{\line(1,1){15}}
\multiput(0,20)(20,0){1}{\circle{6}}
\multiput(0,-20)(20,0){1}{\circle{6}}
\multiput(20,0)(20,0){1}{\circle{6}}
\multiput(80,0)(20,0){1}{\circle{6}}
\multiput(100,0)(20,0){1}{\circle{6}}
\multiput(23,0)(20,0){1}{\line(1,0){14}}
\multiput(63,0)(20,0){1}{\line(1,0){14}}
\multiput(82.85,-1)(0,2){2}{\line(1,0){14.3}}
\multiput(39,0)(4,0){6}{\line(1,0){2}}
\put(90,0){\makebox(0,0){$>$}}
\put(0,-3){\makebox(0,0)[t]{$1$}}
\put(0,37){\makebox(0,0)[t]{$0$}}
\put(20,17){\makebox(0,0)[t]{$2$}}
\put(80,17){\makebox(0,0)[t]{$n\!\! -\!\! 1$}}
\put(100,14){\makebox(0,0)[t]{$n$}}
\put(0,-27){\makebox(0,0)[t]{$n\!\! -\!\! 1$}}
\put(0,10){\makebox(0,0)[t]{$n$}}
\put(24,-7){\makebox(0,0)[t]{$n\!\! -\!\! 2$}}
\put(80,-7){\makebox(0,0)[t]{$1$}}
\put(100,-7){\makebox(0,0)[t]{$0$}}
\put(50,-40){\makebox(0,0)[t]{$B^{(1)}_{n}$}}
\end{picture}
\quad\quad\quad
\begin{picture}(106,85)(-5,-43)
\multiput(2.5,17.5)(0,2){1}{\line(1,-1){15}}
\multiput(2.5,-17.5)(0,2){1}{\line(1,1){15}}
\multiput(0,20)(20,0){1}{\circle{6}}
\multiput(0,-20)(20,0){1}{\circle{6}}
\multiput(20,0)(20,0){1}{\circle{6}}
\multiput(80,0)(20,0){1}{\circle{6}}
\multiput(100,0)(20,0){1}{\circle{6}}
\multiput(23,0)(20,0){1}{\line(1,0){14}}
\multiput(63,0)(20,0){1}{\line(1,0){14}}
\multiput(82.85,-1)(0,2){2}{\line(1,0){14.3}}
\multiput(39,0)(4,0){6}{\line(1,0){2}}
\put(90,0){\makebox(0,0){$<$}}
\put(0,-3){\makebox(0,0)[t]{$1$}}
\put(0,37){\makebox(0,0)[t]{$0$}}
\put(20,17){\makebox(0,0)[t]{$2$}}
\put(80,17){\makebox(0,0)[t]{$n\!\! -\!\! 1$}}
\put(100,14){\makebox(0,0)[t]{$n$}}
\put(0,-27){\makebox(0,0)[t]{$n\!\! -\!\! 1$}}
\put(0,10){\makebox(0,0)[t]{$n$}}
\put(24,-7){\makebox(0,0)[t]{$n\!\! -\!\! 2$}}
\put(80,-7){\makebox(0,0)[t]{$1$}}
\put(100,-7){\makebox(0,0)[t]{$0$}}
\put(50,-40){\makebox(0,0)[t]{$A^{(2)}_{2n-1}$}}
\end{picture}
\end{figure}
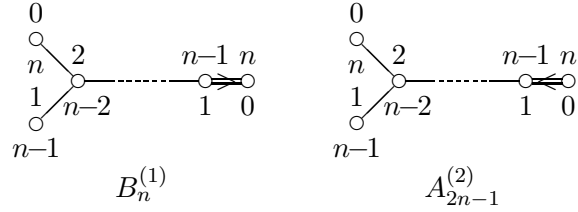

Doing a similar calculation to
Lemmas \ref{lem:w1} and \ref{lem:w2},
we have

\begin{lem}
\label{lem:w3}
Under the bijection $\phi$,
we have the correspondence
($\dot{\Lambda}_0=0$):
\par
Case I.
\begin{align}
\Lambda_a &\mapsto
\begin{cases}
 \dot{\Lambda}_{n-a}-\dot{\Lambda}_{n}
&a=1, n\\
\dot{\Lambda}_{n-a}-
2\dot{\Lambda}_n & a=2,\dots,n-1,\\
\end{cases}\\
\alpha_a &\mapsto
\begin{cases}
 \dot{\alpha}_{n-a}
& a=1,\dots,n-1\\
 -\frac{1}{2}(2\dot{\alpha}_1+\dots
+2\dot{\alpha}_{n-2}+
\dot{\alpha}_{n-1}+\dot{\alpha}_n)& a=n.
\end{cases}
\end{align}
\par
Case II.
\begin{align}
\Lambda_a &\mapsto
\begin{cases}
 \dot{\Lambda}_{n-1}-\dot{\Lambda}_{n}
&a=1\\
\dot{\Lambda}_{n-a}-
2\dot{\Lambda}_n & a=2,\dots,n,\\
\end{cases}\\
\alpha_a &\mapsto
\begin{cases}
 \dot{\alpha}_{n-a}
& a=1,\dots,n-1\\
 -(2\dot{\alpha}_1+\dots
+2\dot{\alpha}_{n-2}+
\dot{\alpha}_{n-1}+\dot{\alpha}_n)& a=n.
\end{cases}
\end{align}
\end{lem}

\begin{lem}
\label{lem:w4}
There is an element $s\in \cW(D_n)$ which acts
on $\dot{\fh}^*$ as follows:
\par
Case I.
\begin{align}
\phi(\Lambda_a)
 &\mapsto 
\begin{cases}
\dot{\Lambda}_a & a=1,\dots,n-2,n\\
\dot{\Lambda}_{n-1}+\dot{\Lambda}_{n} & a=n-1,\\
\end{cases}
\\
\phi(\alpha_{a}) &\mapsto 
\begin{cases}
\dot{\alpha}_a & a=1,\dots,n-1\\
\frac{1}{2}(-\dot{\alpha}_{n-1}+\dot{\alpha}_n) &a=n.\\
\end{cases}
\end{align}
\par
Case II.
\begin{align}
\phi(\Lambda_a)
 &\mapsto 
\begin{cases}
\dot{\Lambda}_a & a=1,\dots,n-2\\
\dot{\Lambda}_{n-1}+\dot{\Lambda}_{n} & a=n-1\\
2\dot{\Lambda}_n & a=n,\\
\end{cases}
\\
\phi(\alpha_{a}) &\mapsto 
\begin{cases}
\dot{\alpha}_a & a=1,\dots,n-1\\
-\dot{\alpha}_{n-1}+\dot{\alpha}_n &a=n.\\
\end{cases}
\end{align}
\end{lem}

Accordingly,
we set
\par
Case I.
\begin{align}
x_a&=e^{\dot{\Lambda}_a}\ (a=1,\dots,n-2,n),\
e^{\dot{\Lambda}_{n-1}+\dot{\Lambda}_{n}}\ (a=n-1),\\
\label{eq:ya1}
y_a&=e^{-\dot{\alpha}_a}\ (a=1,\dots,n-1),\
e^{(\dot{\alpha}_{n-1}-\dot{\alpha}_n)/2}\ (a=n).
\end{align}
\par
Case II.
\begin{align}
x_a&=e^{\dot{\Lambda}_a}\ (a=1,\dots,n-2),\
e^{\dot{\Lambda}_{n-1}+\dot{\Lambda}_{n}}\ (a=n-1),\
e^{2\dot{\Lambda}_{n}}\ (a=n),\\
\label{eq:ya2}
y_a&=e^{-\dot{\alpha}_a}\ (a=1,\dots,n-1),\
e^{\dot{\alpha}_{n-1}-\dot{\alpha}_n}\ (a=n).
\end{align}
Then, the relation (\ref{eq:yx4}) is preserved.
Define the $\dot{\fg}$-characters of $W\am(\zeta)$
in the same way as Definition \ref{def:fc2}.
Then,
the normalized $\fg_0$-character and the
normalized $\dot{\fg}$-character
of $W\am(\zeta)$ coincide
as polynomials of $y$.
Thus, Conjecture \ref{conj:KR2}
for the normalized $\fg_0$-characters
  is applied
for the normalized $\dot{\fg}$-characters as well.
So far, 
the situation is parallel to the
$C_n$ case for $A^{(2)}_{2n}$.
{}From now on,
the situation is parallel to the
$B_n$ case for $A^{(2)}_{2n}$.
The following relations are
easily derived from the explicit expressions
of the Weyl denominators for $B_n$, $C_n$, $D_n$:
\begin{align}
\label{eq:denom12}
\prod_{\alpha\in {\Delta}^{B_n}_+}(1-e^{-\alpha})
&=\prod_{a=1}^n\Bigl(1-\prod_{k=a}^{n} y_k\Bigr)
\prod_{\alpha\in \Delta^{D_n}_+}(1-e^{-\alpha}),\\
\label{eq:denom13}
\prod_{\alpha\in {\Delta}^{C_n}_+}(1-e^{-\alpha})
&=\prod_{a=1}^n\Bigl(1-y_n^{-1}\prod_{k=a}^{n} y_k^2\Bigr)
\prod_{\alpha\in \Delta^{D_n}_+}(1-e^{-\alpha}),
\end{align}
where the equalities 
hold under the following identifications:
(\ref{eq:ya0}) for the LHSs,
(\ref{eq:ya1}) for the RHS of (\ref{eq:denom12}),
(\ref{eq:ya2}) for the RHS of (\ref{eq:denom13})
under the label in Figure \ref{fig:d3}.
We conclude that 
Conjecture \ref{conj:KR2} 
for  $B^{(1)}_n$ and $A^{(2)}_{2n-1}$  is equivalent to

\begin{conj}
\label{conj:KR7}
(i) 
For  $B^{(1)}_{n}$, the formula
\begin{align}
\label{eq:krc7}
\cQ^\nu(y)=
\frac{
\cK^\nu_{D,G}(y)
\prod_{a=1}^n\bigl(1-\prod_{k=a}^{n} y_k\bigr)^{-1}}
{\displaystyle
\prod_{{\alpha} \in {\Delta}^{D_n}_+}(1-e^{-{\alpha}})}
\end{align}
holds
for the $D_n$-characters of the KR-modules,
where $y$ is specified as (\ref{eq:ya1}).
\par
(ii)
For  $A^{(2)}_{2n-1}$,
the formula
\begin{align}
\label{eq:krc8}
\cQ^\nu(y)=
\frac{
\cK^\nu_{D,G}(y)
\prod_{a=1}^n\bigl(1-y_n^{-1}\prod_{k=a}^{n} y_k^2\bigr)^{-1}}
{\displaystyle
\prod_{{\alpha} \in {\Delta}^{D_n}_+}(1-e^{-{\alpha}})}
\end{align}
holds
for the $D_n$-characters of the KR-modules,
where $y$ is specified as (\ref{eq:ya2}).
\end{conj}

The manifest polynomial expressions of the numerators
in the RHSs of (\ref{eq:krc4}), (\ref{eq:krc7}),
and (\ref{eq:krc8})
for $\cQ\am(y)$ are available in \cite{HKOTT}
with some other examples.

\subsection{Related works}
\label{subsec:rel1}

Below we list the related works
on Conjectures  \ref{conj:KR1}
and \ref{conj:KR2}--\ref{conj:KR4}
mostly chronologically.
However, the list is by no means complete.
The series $\cK^\nu_{D,G}(y)$ in
(\ref{eq:krc3}) admits a natural
$q$-analogue called the {\em fermionic formula\/}.
This is another fascinating subject, but we do not cover
it here.
See \cite{BS,HKOTY,HKOTT} and reference therein.
It is convenient to refer the formula
(\ref{eq:krc3})
with the binomial coefficient (\ref{eq:binom1})
as {\em type I},
and the ones with the binomial coefficient
in Remark \ref{rem:binom1}
as {\em type II}.
(In the context of the \XXX-type integrable spin chains,
$N\am$ and $P\am$
represent the numbers of {\em $m$-strings\/} and {\em $m$-holes\/}
of color $a$, respectively.
Therefore one must demand $P\am \ge 0$,
which implies that the relevant formulae are
necessarily of type II.)
The manifest expression of
the decomposition of $\bQ\am$
such as
\begin{align}
\bQ^{(2)}_1= \chi(\Lambda_2)+\chi(\Lambda_5)
\end{align}
is referred as {\em type III\/},
where $\chi(\lambda)$ is the character
of the irreducible $X_n$-module
$V(\lambda)$ with highest weight $\lambda$.
Since there is no  essential distinction 
between these conjectured formulae
for $Y(X_n)$ and $U_q(X^{(1)}_n)$,
we simply refer the both cases  as $X_n$ below.
At this moment, however,
the proofs should be separately given
for nonsimply-laced case \cite{V}.


0 \cite{Be}. Bethe solved the \XXX\/ spin chain of
length $N$ by inventing what is later known
as the Bethe ansatz and the string hypothesis.
As a check of the completeness of his eigenvectors
for the \XXX\/ Hamiltonian,
he proved, in our terminology, the type II formula of $\cQ^\nu(y)$
with $\nu^{(1)}_m=N\delta_{m 1}$ for  $A_1$.
See \cite{F,FT} for a readable exposition
in the framework of the quantum inverse scattering method.

1 \cite{K1,K2}.
Kirillov proposed and proved
the type I formula of the irreducible modules
$V(m\Lambda_a)$ for $A_1$ \cite{K1}
and $A_n$ \cite{K2}. 
The idea of the use of the generating function and the
$Q$-system, which is extended in the present paper,
 originates in this work.

2 \cite{KKR}.
Kerov {\em et al}.\ proposed and proved the type II formula
for $A_n$ by the combinatorial method,
where the bijection between the Littlewood-Richardson
tableaux and the rigged configurations was constructed.

3 \cite{D1}.
Drinfeld claimed that $V(m\Lambda_a)$
can be lifted to a $Y(X_n)$-module,
if the Kac label for $\alpha_a$ in $X^{(1)}_n$ is 1.
These modules are often
called the {\em evaluation modules\/},
and  identified with some KR-modules.
A method of proof is given in \cite{C}
for $U_q(X^{(1)}_n)$.
Therefore, the type III formula
$\bQ\am=\chi(m\Lambda_a)$ holds  for those $a$.
Some of the corresponding $R$-matrices for the
classical algebras,
$X_n=A_n$, $B_n$, $C_n$, $D_n$, 
were obtained earlier in \cite{KRS,R}
by the  {\em reproduction scheme\/}
(also known as the {\em fusion procedure})
in the context of the algebraic Bethe ansatz
method.

4 \cite{OW}.
Ogievetsky and Wiegmann
proposed the type III formula of $\bQ^{(a)}_1$
for some $a$ for the exceptional algebras
{}from the reproduction scheme.

5 \cite{KR,K3}.
Kirillov and Reshetikhin
formulated  the type II
formula for any simple Lie algebra $X_n$.
For that purpose, they vaguely introduced
a family of $Y(X_n)$-modules,
which we identify with the KR modules here.
They proposed  the type II formula for
any $X_n$, and the $Q$-system and
the type III formula for $X_n=B_n$, $C_n$, $D_n$.
The $Q$-system for exceptional algebras $X_n$ was
also proposed  in \cite{K3}.
Due to the long-term absence
of the proofs of the announced results  by
the authors, we regard
these statements as conjectures at our discretion
in this paper.
See Remark \ref{rem:KR2} for further remark.

\begin{rem}
\label{rem:KR2}
Let $X_n=B_n$, $C_n$, $D_n$.
Let $\bQ\am$ and $\cQ\am$ be the $X_n$-character
and the normalized $X_n$-character of the ``KR module" proposed in
\cite{KR}. Then,
one can organize the conjectures in \cite{KR} as follows:
\begin{itemize}
\item[(i)] All  $\bQ\am$'s are given by the type III formula in
\cite{KR}.
\item[(ii)]
The family $(\bQ\am)_{(a,m)\in H}$ satisfies the $Q$-system
(\ref{eq:qsys25}) for $X_n$,
and $\bQ^{(a)}_1$'s $(a=1,\dots,n)$
are given by the type III formula in \cite{KR}.
(Note that the $Q$-system (\ref{eq:qsys25}),
or equivalently (\ref{eq:qsys23}), recursively determines
all  $\bQ\am$'s from the initial data $(\bQ^{(a)}_1)_{a=1}^n$.)
\item[(iii)] Any power $\cQ^\nu$ is  given by the type II formula.
\end{itemize}
As stated in \cite{KR},
one can certainly show
the equivalence between (i) and (ii)
without referring to the KR modules themselves.
See \cite{HKOTY}.
One can also confirm the equivalence between (i) and a
weak version of (iii):
\begin{itemize}
\item[(iii')] All  $\cQ\am$'s are  given by the type II formula.
\end{itemize}
See \cite{Kl} and Appendix A in \cite{HKOTY}.
The family $(\cQ\am)_{(a,m)\in H}$ given by (i)
satisfies the convergence property (\ref{eq:con1}).
Thus, (i), (ii), and (iii') are all equivalent to
\begin{itemize}
\item[(iv)] The family
$(\cQ\am)_{(a,m)\in H}$ is the canonical solution of
the
$Q$-system (\ref{eq:qsys23}).
\end{itemize}
Therefore, as shown in Section \ref{subsec:Rc} (also \cite{KN2}),
they are also equivalent
to
\begin{itemize}
\item[(v)] Any power $\cQ^\nu$ is  given by the type I formula (\ref{eq:krc1}).
\end{itemize}
This is why we call Conjecture \ref{conj:KR1} the Kirillov-Reshetikhin
conjecture.
The equivalence between (iii) and the others has not been proved yet
as we mentioned in Remark \ref{rem:binom1}.
\end{rem}

6 \cite{CP1,CP2}.
Chari and Pressley proved
the type III formula of $\bQ^{(a)}_1$ 
in most cases  for
 $Y(X_n)$ \cite{CP1}, and for $U_q(X^{(1)}_n)$ \cite{CP2},
where the list is complete except for  $E_7$ and  $E_8$.

7 \cite{Ku}. The type III formula of
$\bQ\am$ was proposed for
some $a$ for the exceptional algebras.

8 \cite{Kl}.
Kleber analyzed a combinatorial structure
of the type II formula for the simply-laced algebras.
In particular, it was proved that 
the type III formula of $\bQ\am$
and the corresponding type II formula
are equivalent for $A_n$ and $D_n$.

9 \cite{HKOTY,HKOTT}.
Hatayama {\em et al}.\  gave
a characterization of the type I formula 
as the solution of the $Q$-system 
which are $\bC$-linear combinations of the $X_n$-characters
with the property equivalent to the convergence property
(\ref{eq:con1}).
Using it,
the equivalence of the type III formula of $\bQ\am$
and the  type I formula of $\cQ^\nu(y)$
for the classical algebras was shown \cite{HKOTY}.
In \cite{HKOTT},
the type I and type II formulae, and the $Q$-systems
 for the twisted algebras $U_q(X^{(r)}_N)$
were proposed.
The type III formula of $\bQ\am$ for $A^{(2)}_{2n}$,
$A^{(2)}_{2n-1}$, $D^{(2)}_{n+1}$, $D^{(3)}_4$
was also proposed,
and the equivalence to the type I formula was
shown in a similar way to the untwisted case.

10 \cite{KN1,KN2}.
The second formula in Conjecture \ref{conj:KR3} was
proposed and proved for $A_1$ \cite{KN1} from the
formal completeness of the \XXZ-type Bethe vectors.
The same formula was proposed for $X_n$,
and the equivalence to the type I formula was proved
\cite{KN2}.
The type I formula is formulated in the form
(\ref{eq:fc1}), and the characterization of type I formula 
in \cite{HKOTY} was simplified
as the solution of the $Q$-system
with the convergence property (\ref{eq:con1}).

11 \cite{C}.
Chari proved the type III formula of $\bQ\am$
for $U_q(X^{(1)}_n)$ 
 for any $a$ for the
classical algebras, and for some $a$ for the
exceptional algebras.

12 \cite{OSS}.
Okado {\em et al.\ }constructed bijections
between the rigged configurations and the 
crystals 
(resp.\ virtual crystals)
corresponding to $\cQ^\nu(y)$,
with $\nu\am=0$ for $m>1$, for $C^{(1)}_n$ and $A^{(2)}_{2n}$
(resp.\ $D^{(2)}_{n+1}$).
As a corollary, the type II formula of those
$\cQ^\nu(y)$ was proved for $C^{(1)}_n$ and $A^{(2)}_{2n}$.

\par
Assembling all the above results and the indications
to each other,
let us summarize the current status of the Kirillov-Reshetikhin
conjectures into the following theorem.
Here, we mention the results only for the
quantum affine algebra $U_q(X^{(r)}_{N})$ case.
Also, we exclude the isolated results only valid
for small $m$.

\begin{thm}
(i) Conjecture \ref{conj:KR2}
and the type I formula  of $\cQ^\nu(y)$
are valid for
$A^{(1)}_n$, $B^{(1)}_n$, $C^{(1)}_n$, $D^{(1)}_n$.

\par
(ii) The type II formula of $\cQ^\nu(y)$ is valid for
$A^{(1)}_n$
and valid
for those $\nu$
with $\nu\am=0$ for $m>1$
for $C^{(1)}_n$ and $A^{(2)}_{2n}$.
The type II formula of $\cQ\am(y)$ is valid for
the following $a$ in \cite{C}:
any $a$ for $B^{(1)}_n$, $C^{(1)}_n$, $D^{(1)}_n$;
$a=1,6$ for $E^{(1)}_6$; $a=7$ for $E^{(1)}_7$.

\par
(iii) The type III formula of $\bQ\am$ is valid for 
all $a$ for $A^{(1)}_n$, $B^{(1)}_n$, $C^{(1)}_n$,
$D^{(1)}_n$,
and for those $a$ listed in \cite{C} for 
$E^{(1)}_6$, $E^{(1)}_7$, $E^{(1)}_8$,
$F^{(1)}_4$, $G^{(1)}_2$.
The formula is found in \cite{C}.
\end{thm}

\begin{appendix}
\section{The Denominator formulae}

We give a proof of Proposition \ref{prop:denom11}.
The proof is divided into three steps.

\label{sec:denom}
\subsection{Step 1. Reduction of the denominator formula}
In Steps 1 and 2,
 we consider the unspecialized infinite  $Q$-system
(\ref{eq:qsys10}),
and we assume that $D$ and $G$ satisfy the
condition (KR-II) in Definition \ref{defn:KR1}.

\par
For a given positive integer $L$,
let $H_L=\{1,\dots,n\}\times\{1,
\dots,L\}$ be the finite subset of $H$
 in Section \ref{subsec:spec1}.
With  multivariables $v_L=(v\am)_{(a,m)\in H_L}$,
$w_L=(w\am)_{(a,m)\in H_L}$,
$z_L=(z\am)_{(a,m)\in H_L}$,
we define the bijection $v_L\mapsto w_L$
around $v=w=0$  (cf.\ (\ref{eq:wv1}))
by
\begin{align}
\label{eq:wv5}
w\am(v_L)&=v\am
\prod_{(b,k)\in H_L}(1-v\bk)^{-G'_{am,bk}},
\end{align}
and the bijection $v_L\mapsto z_L$
around $v=z=0$ by
\begin{align}
\label{eq:zv2}
z\am(v_L)&=w\am(v_L)
\prod_{(b,k)\in H_L}(1-v\bk)^{g_{ab}m}
\\
\label{eq:zv1}
&=v\am\prod_{(b,k)\in H_L}(1-v\bk)^{-G'_{am,bk}+g_{ab}m},
\end{align}
where $g_{ab}$ is the one in (KR-II).
Let us factorize the bijection $w_L\mapsto v_L$
as $w_L\mapsto z_L\mapsto v_L$.
The map $w_L\mapsto z_L$ is described as
\begin{align}
\label{eq:zw5}
z\am(w_L)&=w\am
\prod_{b=1}^n
(Q_b(w_L))^{-g_{ab}m},
\quad
Q_b(w_L):=
\prod_{k=1}^L(1-v\bk(w_L))^{-1}.
\end{align}
By the assumption (KR-II)
and the expression (\ref{eq:zv1}),
 the map $v_L\mapsto z_L$ is
lower-triangular in the sense of 
Example \ref{exmp:tri1}.
Therefore, the following equality holds:
\begin{align}
\label{eq:dwv1}
\det_{H_L}\Bigl(
\frac{w\bk}{v\am}
\frac{\partial v\am}{\partial w\bk}(w_L)
\Bigr)
=
\det_{H_L}\Bigl(
\frac{w\bk}{z\am}
\frac{\partial z\am}{\partial w\bk}(w_L)
\Bigr),
\end{align}
where $\det_{H_L}$ means the abbreviation of
$\det_{(a,m),(b,k)\in H_L}$.
\par
We now simultaneously
specialize the variables $w_L$ and $z_L$
with the variables
$y=(y_a)_{a=1}^n$ and $u=(u_a)_{a=1}^n$ as
(cf.\ (\ref{eq:wy3}))
\begin{align}
\label{eq:zu1}
w\am=w\am(y)=(y_a)^m,
\quad z\am=z\am(u)=(u_a)^m.
\end{align}
This specialization is compatible with (\ref{eq:zw5})
and the map $y\mapsto u$,
\begin{align}
\label{eq:uy1}
u_a(y)=y_a \prod_{b=1}^n (q_{b}(y))^{-g_{ab}},
\quad
q_{b}(y):=Q_b(w_L(y)).
\end{align}

\begin{prop}
Let
$G'_L=(G'_{am,bk})_{(a,m),(b,k)\in H_L}$
be the $H_L$-truncation of $G'$,
$K^0_{I_L,G'_L}(w_L)$ be the one
in (\ref{eq:pse2}),
and $\cK^0_{I_L,G'_L}(y):=K^0_{I_L,G'_L}
(w_L(y))$ be its specialization by (\ref{eq:zu1}).
Then, the formula (\ref{eq:pse2}) reduces
to
\label{prop:k03}
\begin{align}
\label{eq:k06}
\cK^0_{I_L,G'_L}(y)&=
\det_{1\leq a,b\leq n}\Bigl(
\frac{y_b}{u_a}
\frac{\partial u_a}{\partial y_b}
(y)
\Bigr)
\prod_{a=1}^n q_{a}(y).
\end{align}
\end{prop}
\begin{proof}
Because of (\ref{eq:dwv1}),
it is enough to prove the equality
\begin{align}
\label{eq:dzw1}
\det_{H_L}\Bigl(
\frac{w\bk}{z\am}
\frac{\partial z\am}{\partial w\bk}(w_L(y))
\Bigr)
&=
\det_{1\leq  a,b\leq n}\Bigl(
\frac{y_b}{u_a}
\frac{\partial u_a}{\partial y_b}(y)
\Bigr).
\end{align}
We remark that
\begin{align}
\label{eq:dx2}
y_a\frac{\partial}{\partial y_a}
&=
\sum_{m=1}^L m w\am
\frac{\partial}{\partial w\am},\\
\label{eq:detf3}
\det_{H_L}(\delta_{am,bk}+m\alpha_{abk})
&=
\det_{1\leq a,b\leq n}
\Bigl(\delta_{ab}+\sum_{k=1}^L k\alpha_{abk}
\Bigr),
\end{align}
where $\alpha_{abk}$ are arbitrary constants depending
on $a$, $b$, $k$.
Set $F_a(w_L)=
\prod_{b=1}^n (Q_b(w_L))^{-g_{ab}}$.
Then, (\ref{eq:dzw1}) is obtained as
\begin{align*}
\text{(LHS)}
&=
\det_{H_L}\Bigl(
\delta_{am,bk}+mw\bk
\frac{\partial }{\partial w\bk}
\log F_a(w_L(y))
 \Bigr)\\
&=
\det_{1\leq a,b\leq n}
\Bigl(
\delta_{ab}
+\sum_{k=1}^L k w\bk
\frac{\partial }{\partial w\bk}
\log F_a(w_L(y))
\Bigr)\\
&=
\det_{1\leq a,b\leq n}
\Bigl(
\delta_{ab}
+y_b\frac{\partial }{\partial y_b}
\log F_a(w_L(y))
\Bigr)\\
&=
\det_{1\leq  a,b\leq n}
\Bigl(
\frac{y_b}{u_a}
\frac{\partial u_a}{\partial y_b}(y)
\Bigr),
\end{align*}
where we used (\ref{eq:zw5}),
(\ref{eq:detf3}),
(\ref{eq:dx2}), and (\ref{eq:uy1}).
\end{proof}

\subsection{Step 2. Change of variables}
We introduce the change of
the variables $y$ and $u$
in (\ref{eq:zu1}) to $x=(x_a)_{a=1}^n$
 and $\bq=(\bq_a)_{a=1}^n$ as
\begin{align}
\label{eq:yx1}
y_a(x)&=\prod_{b=1}^n (x_b)^{-g_{ab}},
\quad
u_a(\bq)=\prod_{b=1}^n (\bq_b)^{-g_{ab}}.
\end{align}
Thus, if $f(y)$ is a power series of $y$, then $f(y(x))$ is
a Laurent series of $x$
because of the assumption in (KR-II)
 that $g_{ab}$'s are integers.
This specialization is compatible with (\ref{eq:uy1})
and the map $x\mapsto \bq$,
\begin{align}
\label{eq:sx1}
\bq_a(x)&=x_a q_{a}(y(x)).
\end{align}
Let us summarize all the maps and variables
 in a diagram:
\begin{align}
\begin{matrix}
v & 
\overset{(\ref{eq:zv2})}{\longleftrightarrow}
 & z & 
\overset{(\ref{eq:zw5})}{\longleftrightarrow}
 & w\\
& {\scriptstyle (\ref{eq:zu1})} \hskip-15pt&\uparrow
& &\uparrow & \hskip-8pt{\scriptstyle (\ref{eq:zu1})}\\
&&u&
\overset{(\ref{eq:uy1})}{\longleftrightarrow}
&y\\
&{\scriptstyle (\ref{eq:yx1})}\hskip-15pt &\uparrow
& &\uparrow& \hskip-8pt{\scriptstyle (\ref{eq:yx1})}\\
&&\bq&
\overset{(\ref{eq:sx1})}{\longleftrightarrow}
&x
\end{matrix}
\end{align}

With these changes of variables,
(\ref{eq:k06}) becomes the Jacobian of $\bq(x)$:
\begin{prop}
Let $\cK^0_{I_L,G'_L}(y)$ be the one in
Proposition \ref{prop:k03},
and let $\bK^0_{I_L,G'_L}(x)
:=\cK^0_{I_L,G'_L}(y(x))$.
Then, the formula
\label{prop:k01}
\begin{align}
\label{eq:k05}
\bK^0_{I_L,G'_L}(x)&=
\det_{1\leq a,b\leq n}
\Bigl(\frac{\partial \bq_a}
{\partial x_b}(x)\Bigr)
\end{align}
holds.
\end{prop}

\begin{proof}
By (\ref{eq:yx1}), we have
\begin{align}
\label{eq:dsz1}
\det_{1\leq a,b\leq n}\Bigl(
\frac{\bq_b}{ u_a}
\frac{\partial u_a}{\partial \bq_b}
\Bigr)
&=
\det_{1\leq a,b\leq n}\Bigl(
\frac{x_b}{y_a}
\frac{\partial y_a}{\partial x_b}
\Bigr)
=
\det_{1\leq a,b\leq n}
(-g_{ab})\neq 0.
\end{align}
Using Proposition \ref{prop:k03},
(\ref{eq:sx1}), and (\ref{eq:dsz1}),
we obtain
\begin{align}
\begin{split}
\bK^0_{I_L,G'_L}(x)
&=
\det_{1\leq a,b\leq n}\Bigl(
\frac{y_b}{u_a}
\frac{\partial u_a}{\partial y_b}(y(x))
\Bigr)
\prod_{a=1}^n q_a(y(x))\\
&=
\det_{1\leq a,b\leq n}\Bigl(
\frac{x_b}{\bq_a}
\frac{\partial \bq_a}{\partial x_b}(x)
\Bigr)
\prod_{a=1}^n q_a(y(x))
=
\det_{1\leq a,b\leq n}
\Bigl(
\frac{\partial \bq_a}{\partial x_b}(x)
\Bigr).
\end{split}
\end{align}
\end{proof}

\subsection{Step 3. Denominator formula
for the $Q$-systems for KR type}

Now we are ready to prove Proposition \ref{prop:denom11};
namely,

\begin{prop}
\label{prop:denom12}
Let $\bK^0_{D,G}(x):=
\cK^0_{D,G}(y(x))$, where
$\cK^0_{D,G}(y)$ is
the denominator in (\ref{eq:qkr8})
for the $Q$-system of KR type (\ref{eq:qsys14}).
Then, the formula
\begin{align}
\label{eq:hk1}
\bK^0_{D,G}(x)=
\det_{1\leq a,b\leq n}
\Bigl(
\frac{\partial \bQ\sa_1}{\partial x_b}(x)
\Bigr)
\end{align}
holds, where we set $\bQ\sa_1(x):=x_a\cQ^{(a)}_1(y(x))$
for the canonical solution $(\cQ\am(y))_{(a,m)\in H}$
of  (\ref{eq:qsys14}).
\end{prop}

\begin{proof}
We recall the following four facts:

\par
Fact 1: By (\ref{eq:kl1}) and (\ref{eq:wycom1}), we have
\begin{align}
\label{eq:k010}
\cK^0_{D,G}(y)
\equiv
\cK^0_{I_L,G'_L}(y)
\mod J_L.
\end{align}

\par
Fact 2:
By Theorem \ref{thm:cs4} and the proof therein,
the canonical solution 
$(\cQ\am(y))_{(a,m)\in H}$ of (\ref{eq:qsys14})
and the solution $(\cQ^{\prime (a)}_m(y))_{(a,m)\in H}$
of the corresponding standard $Q$-system are related as
\begin{align}
\cQ^{\prime(a)}_m(y)=
\frac{\cQ\amp(y)
\cQ\amm(y)}{(\cQ\am(y))^2}.
\end{align}

\par
Fact 3:
By Propositions \ref{prop:qsys1}, \ref{prop:iq1},
 and (\ref{eq:wycom1}),
the series $q_b(y)$ in (\ref{eq:uy1})
satisfies
\begin{align}
q_a(y)\equiv\prod_{m=1}^L (\cQ^{\prime(a)}_m(y))^{-1}
 \mod J_L,
\end{align}
where $\cQ^{\prime(a)}_m(y)$ is the one in 
Fact 2.
Note that $q_b(y)$ depends on $L$.

\par
Fact 4:
By the proof of Proposition \ref{prop:conv1},
it holds that
\begin{align}
\cQ\aL(y)\equiv \cQ\aLp(y)
\mod J_L.
\end{align}
Combining Facts 2--4, we immediately have
$q_a(y)\equiv \cQ^{(a)}_1(y)$ mod $J_L$.
Thus, $\lim_{L\to \infty}q_a(y)=\cQ^{(a)}_1(y)$ holds.
Therefore, taking the limit $L\rightarrow \infty$
of (\ref{eq:k06}) with the help of Fact 1, we obtain
\begin{align}
\label{eq:ku1}
\cK^0_{D,G}(y)=&\,
\det_{1\leq a,b\leq n}\Bigl(
\frac{y_b}{U_a}
\frac{\partial U_a}{\partial y_b}
(y)
\Bigr)
\prod_{a=1}^n \cQ^{(a)}_{1}(y),\\
\label{eq:ku2}
U_a(y):=&\, y_a\prod_{b=1}^n
(\cQ^{(b)}_1(y))^{-g_{ab}}.
\end{align}
The equality (\ref{eq:hk1}) is obtained  from
(\ref{eq:ku1}) in the same way
 as the proof of Proposition \ref{prop:k01}.
\end{proof}

\end{appendix}

\end{document}